\documentclass[AMA,STIX1COL]{WileyNJD-v2}

\articletype{Research Article}%

\received{14th September 2018}
\revised{22nd February 2019}
\accepted{18th March 2019, International Journal for Numerical Methods in 
Engineering}

\usepackage[utf8]{inputenc}
\usepackage[T1]{fontenc}
\usepackage{amsmath}
\usepackage{amssymb}
\usepackage{amsfonts}
\usepackage{soul}
\usepackage{color}
\usepackage{units}
\usepackage{subcaption}
\usepackage{graphicx}
\usepackage{url}
\usepackage{textcomp}
\usepackage{xspace}

\renewcommand{\P}{\mathcal{P}}
\renewcommand{\L}{\mathcal{L}}
\newcommand{\R}{\mathbb{R}}
\newcommand{\T}{\mathsf{T}}
\newcommand{\eh}{\hat{e}}
\renewcommand{\b}{\boldsymbol}
\newcommand{\w}{\b{w}}
\newcommand{\eps}{\varepsilon}
\renewcommand{\phi}{\varphi}
\newcommand{\dpar}[2]{\frac{\partial #1}{\partial #2}}
\newcommand{\Ns}{\mathcal{N}}
\newcommand{\dr}{\delta r}
\newcommand{\drI}{\dr^{(0)}}
\newcommand{\drup}{\dr^u}
\newcommand{\drlo}{\dr^\ell}

\newcommand{\prf}{\emph{p}-refinement\xspace}
\newcommand{\hrf}{\emph{h}-refinement\xspace}
\newcommand{\rrf}{\emph{r}-refinement\xspace}
\newcommand{\hrrf}{\emph{h,r}-refinement\xspace}

\usepackage{algpseudocode}  
\usepackage{algorithm}      
\algnewcommand\algorithmicto{\textbf{to}}
\algnewcommand\algorithmicin{\textbf{in}}
\algnewcommand\algorithmicforeach{\textbf{for each}}
\algrenewtext{For}[3]{\algorithmicfor\ #1 $\gets$ #2\ \algorithmicto\ #3\ \algorithmicdo}
\algdef{S}[FOR]{ForEach}[2]{\algorithmicforeach\ #1\ \algorithmicin\ #2\ \algorithmicdo}

\begin{document}

\title{Adaptive Radial Basis Function-generated Finite Differences method for
contact problems}

\author[1,2]{Jure Slak}
\author[1]{Gregor Kosec*}

\address[1]{\orgdiv{Parallel and Distributed Systems Laboratory}, \orgname{``Jožef Stefan'' Institute}, \orgaddress{\country{Slovenia}}}
\address[2]{\orgdiv{Faculty of Mathematics and Physics}, \orgname{University of Ljubljana}, \orgaddress{\country{Slovenia}}}

\corres{*Gregor Kosec, \email{gregor.kosec@ijs.si}}

\presentaddress{
Department E6 \\
``Jožef Stefan'' Institute \\
Jamova ulica 39 \\
1000 Ljubljana \\
Slovenia
}

\abstract[Summary]{
  This paper proposes an original adaptive refinement framework using Radial
  Basis Functions-generated Finite Differences method.
  Node distributions are generated with a Poisson Disk Sampling-based
  algorithm from a given continuous density function, which is altered during
  the refinement process based on the error indicator.
  All elements of the proposed adaptive strategy rely only on
  meshless concepts, which leads to great
  flexibility and generality of the solution procedure.
  The proposed framework
  is tested on four gradually more complex contact problems, governed by the
  Cauchy-Navier equations.  First, a disk under pressure is considered and the
  computed stress field is compared to the closed form solution of the problem
  to assess the basic behaviour of the algorithm and the influence of free
  parameters.  Second, a Hertzian contact problem, also with known closed form
  solution, is studied to analyse the proposed algorithm with an ad-hoc error
  indicator and to test both refinement and derefinement. A
  contact problem, typical for fretting fatigue, with no known closed form
  solution is considered and solved next. It is demonstrated that the proposed
  methodology can be used in practical application and produces results
  comparable with FEM without the need for manual refinement or any human
  intervention. In the last case, generality of the proposed approach is
  demonstrated by solving a \mbox{3-D} Boussinesq's problem of the concentrated
  normal traction acting on an isotropic half-space.
  } 


\maketitle

\section{Introduction}
Solving partial differential equations (PDEs) with Radial basis functions (RBFs)
began in 1990 when Kansa suggested that besides scattered data interpolation,
RBFs could also be used for solving PDEs~\cite{kansa1990multiquadrics}. In the
beginning, only global collocation methods were used, achieving high order
convergence and spectral accuracy~\cite{franke1998solving}, however, at the cost
of high computational complexity and possible ill-conditioning in solution of
the global system, especially as the number of computational nodes increases.
Shortcomings of the global methods encouraged the development of local
approaches and soon several meshless methods based on approximation over local
overlapping support domains appeared~\cite{sarler, prax}. A popular variant
among many local meshless methods is the Radial Basis Function-generated Finite
Differences (RBF-FD) method~\cite{tolstykh2003using}, which uses finite
difference-like collocation weights on an unstructured set of nodes. The method
has been successfully used in several problems and is still actively
researched~\cite{Fornberg_Flyer_2015, bayona2017role, slak2018refined}.In the
field of solid mechanics, where problems are traditionally tackled with the
Finite Element Method (FEM)~\cite{fem_solid}, meshless methods surfaced as a
response to the cumbersome meshing of realistic 3D domains required by
FEM~\cite{slak2018refined, meshless_solid, mavric, Liu}. The development began
with weak form methods such as the Element Free Galerkin
Method~\cite{belytschko1994element} and Meshless Local Petrov Galerkin
Method~\cite{atluri1998new}, followed by strong form collocation approach that
emerged with Finite Point Method (FPM)~\cite{onate2001finite} and further
developed with similar local strong form methods~\cite{slak2018refined, mavric}.
One of the most attractive features of meshless methods is that they do not
require an underlying mesh and fully define relations between nodes directly
through their relative inter nodal positions, which allows to reduce the
problem of constructing a mesh into a much simpler problem of positioning
scattered nodes~\cite{Kos1,fornberg2015fast}. This simplification makes
meshless methods good candidates for adaptive numerical analysis.

Adaptive refinement of domain discretisation is indispensable in problems where
varying precision and resolution of the numerical solution are required in
different areas of the computational domain. A typical example of such a problem
in linear elasticity is a contact of two bodies where extremely high stress
concentrations are present under and around the contact area. Solving such
problems on uniform grids is intractable due to the amount of time and
computational resources required. To obtain a numerical solution with
satisfactory precision in the contact area within a reasonable time frame,
significantly more nodes have to be placed under the contact than elsewhere.
Since the areas of high stresses are generally not known a priori, and to
eliminate the need for human intervention in construction of discretisations, a
commonly adopted approach in such situations is to refine the discretisation
iteratively, gradually increasing the resolution of the numerical solution in
areas of interest with the ultimate goal to ensure relatively uniform error
distribution throughout the domain, or in many cases solely to solve the
problem. In addition to the PDE solution procedure itself, two more
modules are required in an adaptive method, namely an indicator of the solution
quality (an error indicator), and an appropriate refinement strategy.

The most widely accepted error indicator has been introduced by Zienkiewicz and
Zhu in 1987~\cite{Zienkiewicz_1} for FEM solutions of elasticity problems. The
basic idea behind the proposed indicator, also known as ZZ or Z$^2$ is based on
the assumption that the error of the solution can be approximated as the
difference between the numerically obtained solution and a recovered
solution,~i.e.~a more accurate solution computed by appropriate
post-processing. This approach has been later successfully extended to the
meshless solutions of elasticity problems, both weak form,
using meshless finite volume method~\cite{Ebrahimnejad}, NS-PIM~\cite{Tang1},
and strong form, using FPM~\cite{Angulo}. In RBF-FD
context a ZZ type of error indicator has been discussed in a solution of
Laplace equation in~\cite{Oanh}. An alternative class of error indicators
is available for commonly used least squares-based meshless methods,
which relies on the least squares approximation residual~\cite{Sang}.
The residual based error indicator has been successfully used
in a meshless solution of an elasticity problem with a Discrete Least Squares
Meshless method~\cite{afshar2011node}. Furthermore, for specific problems,
ad-hoc meshless error indicators have been reported in literature,
inspired by a physical interpretation of the solution. These include indicators
such as a maximal difference of values within the support domain, effectively
acting as a maximal first derivative of the considered field~\cite{Davydov}, or
variance of the field values within the support domain~\cite{kos2}.

Once the knowledge about the error of the solution is extracted, two
conceptually different adaptive approaches can be employed, namely
the \prf or \hrrf.
In the \prf scheme solution accuracy is varied by changing the order of
the approximation. In meshless methods this results in an increase of the
number of shape functions, modifications of the shape functions or an increase
of the influence domain~\cite{Stevens}. Second approach is to modify
discretization of the domain, either with the \rrf scheme, where nodes are
repositioned to improve the approximation, leaving the total number of the nodes
unchanged, or with \hrf, which influences the total number of nodes.
Successful \rrf has been demonstrated in a meshless solution of  elasticity
problems with discrete least square method~\cite{afshar2011node},
on the solution of an infinite plate with a circular hole and a circular disk
under pressure. On the other hand, \hrf changes the total number of nodes by
placing more points in areas where
needed and removing them from areas that are overpopulated or by recreating a
new distribution of nodes based on the information from the original
distribution. In mesh based methods, the adding and removing approach is also
known as \emph{mesh enrichment}, while creating new distribution approach is
commonly referred to as a \emph{re-meshing}~\cite{Ebrahimnejad}.
In the meshless context, \hrf has been successfully used with the global
RBF Collocation Method~\cite{Libre} in a solution of nearly singular PDEs,
as well as with local strong form meshless methods in a solution of Burgers'
equation~\cite{kos2} and torsion problem~\cite{Liu1}. It has also been
successfully applied in meshless solutions of elasticity
problems~\cite{slak2018refined, Ebrahimnejad, Tang1, Afshar, Raz1}.

Even though the existence of refinement algorithms seems to be an inherent
property of all meshless methods, one has to be careful when interfering with
the nodal configuration to avoid potential negative effects such as ill
conditioning and other possible complications that arise with non smooth nodal
distributions~\cite{driscoll2007adaptive,fornberg2007runge}. This task becomes
even more challenging in a fully adaptive approach where the refinement
algorithm and numerical solution have to be robust enough to fulfil the demands
of the error indicator. In several published \hrf meshless solutions, the
authors used Voronoi diagrams, the dual construction to Delaunay
triangulations, for insertion of new nodes~\cite{Sang, Ebrahimnejad, Tang1,
Afshar, Angulo}. This approach is often not feasible for 3-D problems
and can be computationally demanding.

Davydov and Oanh~\cite{Davydov} discussed a pure meshless \hrf algorithm
for strong form RBF-FD solution of Poisson equation. In their proposed
algorithm, new nodes are added at midpoints between existing nodes.
Since the algorithm does not guarantee smooth distributions, required by
RBF-FD~\cite{fornberg2015fast}, authors had to use special stencil selection
algorithms to improve stability~\cite{Oanh}. Similar \hrf algorithm
achieving a ratio of $2^{17}$
between the densest and coarsest parts of the discretisation
has been recently applied to elasticity
problems~\cite{slak2018refined}. No special stencil selection algorithms were
needed in this case; however, a few steps of regularisation were needed
to sufficiently improve the distribution.


This paper continues the evolution of adaptive meshless methods in the field of
linear elasticity by presenting a robust purely meshless adaptive RBF-FD
solution of contact problems in absence of meshing at any point of the solution
procedure. The refinement procedure is also independent on the dimensionality
of the considered domain, which is demonstrated by solving numerical examples
from linear elasticity with pronounced differences in stress within the domain
in 2-D and 3-D. The proposed methodology uses \emph{re-meshing} type \hrf with
an ad-hoc error indicator.

The performance of the proposed adaptive method is demonstrated on four cases,
governed by Navier-Cauchy equations, gradually increasing in complexity. In the
first case, a disk under pressure is considered, which can be solved in closed
form and is therefore ideal for testing, since the exact error can be computed
during refinement. The second case considers Hertzian contact where an ad-hoc
error indicator type~\cite{Davydov,kos2} based on the characteristics of the
underlying problem is used. The third case deals with
numerical computation of stresses under fretting fatigue conditions, i.e.\ a
contact with slip-stick zone that has been recently investigated with finite
element analysis~\cite{pereira2016on}. The last case serves as a presentation of
the generality and effectiveness of the proposed algorithm in 3-D by solving
the Boussinesq's problem of concentrated normal traction acting on an
isotropic half-space.

All examples in this paper were computed using the in-house open source Medusa
library~\cite{medusa}. All research
data is freely available in a git repository~\cite{git}.

The rest of the paper is organized as follows: the RBF-FD method is described in
section~\ref{sec:rbffd}, the adaptivity algorithm and other relevant components
are described in section~\ref{sec:adaptivity-main}, the results on four
increasingly difficult cases are presented in section~\ref{sec:results} and
finally, the conclusions are presented in section~\ref{sec:conclusions}.

\section{RBF-generated finite differences}
\label{sec:rbffd}
The RBF-FD method is used to obtain the discrete formulation of a PDE.
Consider an elliptic boundary value problem with Dirichlet boundary conditions
\begin{align}
\L u = f \quad &\text{ on } \Omega  \\
u = u_0 \quad &\text{ on } \partial \Omega,
\end{align}
where $f$ and $u_0$ are known functions.
To obtain a discrete representation, $N$ nodes are placed in $\Omega$,
of which $N_i$ nodes in the interior and $N_b$ nodes on the boundary.
Let the set of all discretisation nodes be denoted by $\P$.
Each node $p_i$ is assigned $n_i$ neighbours (including itself), denoted $\Ns(p_i)$.
These neighbours are commonly called the node's \emph{support domain}
or, analogous to the Finite Difference Method (FDM), its \emph{stencil}.

\begin{figure}[h!]
  \centering
  \includegraphics[width=0.4\linewidth]{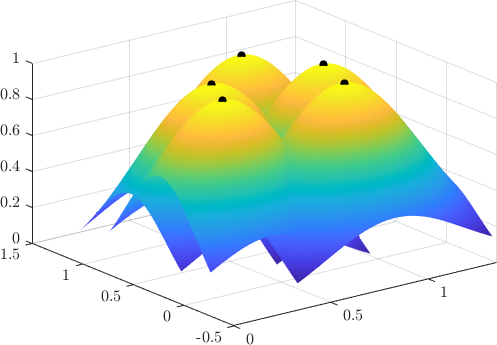}
  \caption{Gaussian radial basis functions with $\sigma_b = 0.5$ over a set of randomly generated nodes.}
  \label{fig:rbfs}
\end{figure}

The operator $\L$ at a point $p_i$ is approximated as a weighted linear
combination of function values at stencil points
\begin{equation}
(\L u)(p_i) \approx \sum_{p_j \in \Ns(p_i)} w^i_j u(p_j).
\label{eq:approx}
\end{equation}
To determine the weights $w^i_j$, exactness of~\eqref{eq:approx} is imposed for
a certain set of functions. A common choice are monomials, as used in the FDM
and FPM, however, in RBF-FD method, as the name suggests, the Radial Basis
Functions (RBFs) are used. RBFs, positioned in centres $P = \{p_1, \ldots, p_n\}$ are functions
\begin{equation}
\{ \phi_i := \phi(\|\cdot - p_i\|), \text{ for } p_i \in P \},
\end{equation}
where $\phi\colon[0, \infty) \to \R$ is an arbitrary function.
In this work Gaussian RBFs will be used, defined as
\[ \phi_i(p) = \exp(-(\|p-p_i\|/\sigma_b)^2), \]
where $\sigma_b$ stands for shape parameter. The Gaussian RBF are illustrated in Fig.~\ref{fig:rbfs}.
Many other types are known and used, with various recommendations~\cite{schaback2001characterization}.

Imposing exactness of~\eqref{eq:approx} for $\phi_k$, we can write
\begin{equation}
(\L \phi_k)(p_i) = \sum_{p_j \in \Ns(p_i)} w^i_j \phi_k(p_j),
\label{eq:exactness}
\end{equation}
for all $p_k \in \Ns(p_i)$. Assembling all equations given by~\eqref{eq:exactness} in
a matrix, the following system of equations is obtained:
\begin{equation}
\begin{bmatrix}
\phi(\|p_{j_1} - p_{j_1}\|) & \cdots & \phi(\|p_{j_{n_i}} - p_{j_1}\|) \\
\vdots & \ddots & \vdots \\
\phi(\|p_{j_1} - p_{j_{n_i}}\|) & \cdots & \phi(\|p_{j_{n_i}} - p_{j_{n_i}}\|) \\
\end{bmatrix}
\begin{bmatrix}
w^i_{j_1} \\ \vdots \\ w^i_{j_{n_i}}
\end{bmatrix}
=
\begin{bmatrix}
(\L \phi_{j_1})(p_i) \\
\vdots \\
(\L \phi_{j_{n_i}})(p_i) \\
\end{bmatrix},
\label{eq:matrix}
\end{equation}
where $j_k$ are indices of nodes in $\Ns(p_i)$.  Equation~\eqref{eq:matrix} is a
system of $n_i$ linear equations, compactly written as $A_i\w_i = b_i$, where
matrix $A_i$ is symmetric.  Furthermore, the properties of such collocation
matrices arising from radial basis functions have been well
studied~\cite{buhmann2003radial,schaback2001characterization}.  If positive
definite RBFs, such as Gaussian RBFs, are used, the matrices $A_i$ are positive
definite as soon as all nodes $p_i$ are distinct~\cite{schaback2001characterization}.
Furthermore, to ensure consistency up to a certain order, RBF-FD approximation can be
augmented with monomials. These consistency constraints are enforced using
Lagrangian multipliers, resulting in a system similar to~\eqref{eq:matrix}.
For example, if a constant 1 is added to the basis set, the system for
computing the weights reads
\begin{equation}
\begin{bmatrix}
\phi(\|p_{j_1} - p_{j_1}\|) & \cdots & \phi(\|p_{j_{n_i}} - p_{j_1}\|) & 1\\
\vdots & \ddots & \vdots & \vdots \\
\phi(\|p_{j_1} - p_{j_{n_i}}\|) & \cdots & \phi(\|p_{j_{n_i}} - p_{j_{n_i}}\|) & 1\\
1 & \cdots & 1  & 0
\end{bmatrix}
\begin{bmatrix}
w^i_{j_1} \\ \vdots \\ w^i_{j_{n_i}} \\ \lambda^i_1
\end{bmatrix}
=
\begin{bmatrix}
(\L \phi_{j_1})(p_i) \\
\vdots \\
(\L \phi_{j_{n_i}})(p_i) \\
0
\end{bmatrix},
\label{eq:matrix-mon}
\end{equation}
effectively adding an additional constraint $\sum_{i=1}^{n_i} w_{j_i} = 0$. The
effects of adding polynomial constraints to RBF-FD approximations have been
studied recently by Bayona et~al.~\cite{bayona2017role}. In cases discussed in
this paper, additional consistency constraints were not needed to obtain
satisfying results. Another possible modification to the RBF-FD above procedure is to include
more stencil points then basis functions, putting basis functions only on the
closest $m_i$ nodes. This turns~\eqref{eq:matrix} into an underdetermined system,
however, it can be solved uniquely by imposing additional condition
of minimizing the norm $\|\w_i\|$ of the weights. The benefit of this approach is that
is less sensitive to perturbations of nodal positions~\cite{slak2018refined}.
This modification will be used in this paper.

After computing the coefficients $\w_i$ for all nodes, a sparse $N\times N$
matrix with $\sum_{i=0}^N n_i$ elements representing the governing problem can
be assembled. The right hand side is computed from given $f$ and $u_0$ and the
sparse system is solved to obtain an approximation for $u$.  Vector PDEs are
treated as a coupled system of scalar PDEs, resulting in a proportionally larger
system.  Boundary conditions that involve differential operators, such as
Neumann or traction boundary conditions, are discretised analogously to operator $\L$ above.

\subsection{Remarks on weight computation and shape parameter $\sigma_b$}
\label{sec:rbffd-shape-comp}

It is worth noting that the selection of the shape parameter $\sigma_b$ requires some attention.
A known problem with computing the weights $\w_i$ is that as the shape parameter
$\sigma_b$ tends towards $\infty$, the accuracy increases, but the matrices $A_i$
become more and more ill conditioned. A common solution is to scale the shape parameter $\sigma_b$
proportionally to the internodal distance and thus keep condition numbers of
$A_i$ constant even as the node density increases. Another possible complication arises from the
Runge Phenomenon, which comes into play with large $\sigma_b$ and variable node densities.
To avoid that, once again, a spatially variable shape parameter has been
used as proposed by Fornberg and Zuev~\cite{fornberg2007runge} and Flyer and
Lehto~\cite{flyer2010rotational}.

A downside to using a spatially variable shape parameter, is that the matrices $A_i$ (Eq.~\ref{eq:matrix})
lose their symmetry properties. To retain the symmetry, the shape parameter is kept constant for
all entries in each $A_i$, but varies with $i$ as
$\sigma_{b,i} = \sigma_b \dr_i$, where $\dr_i$ is the distance to the closest neighbouring node to $p_i$.
Another downside is that scaled shape parameters may cause stagnation errors and failure of
convergence~\cite{flyer2016role}, however, no such problems were detected in
cases discussed in this paper
and scaling the shape parameter sufficed to suppress the Runge Phenomenon and mitigate ill-conditioning.


\section{Adaptivity}
\label{sec:adaptivity-main}

The concept of the adaptively refined solution of a PDE in this work is as
follows: initially discretise the domain with (usually) uniform nodal
distribution, compute approximate solution $u$, compute new desired nodal
density function based on the previous density and the input from the chosen
error indicator, redistribute the nodes and repeat as long as the error is not
within desired range. The refinement procedure presented in this paper is built
on the following preliminaries.
\begin{enumerate}[(i)]
  \item An a-posteriori error indicator $\eh$ is known, where $\eh(p)$ represents
  an indicator of error between an actual and computed value of some desired
  quantity at a point $p$. For convenience we define $\eh(p_i) = \eh_i$ for
  all $p_i \in \P$.  In a more general setting $\eh$ can also represent any
  kind of refinement indicator, with larger values indicating areas where
  refinement is needed and lower values indicating areas where it is
  unnecessary or even derefinement is possible.
  \item Domain is initially filled with a known density $\drI$.
  \item Refinement is an iterative process that alters current radius function
  $\dr$ with respect to the error indicator.
  \item The coarseness limit for $\dr$ is $\drup$, i.e.\ the nodal radius during
  derefinement is bounded by $\drup_i$ from above. \label{enum:deref-bnd}
  \item The nodal radius function $\dr$ is changed with respect to the density
  increase factor ($f_i$) that is evaluated based on the error indicator
  $\hat{e}_i$, the maximal ($\eps$) and minimal ($\eta$) allowed error, and
  maximal ($M$) and minimal ($m$) value of the error indicator. The increase
  factor is limited to $f\in[1/\beta, \alpha]$, where $\alpha, \beta \geq 1$.
\end{enumerate}

In the following subsections the parts of the above procedure are described in more detail.

\subsection{Node placement}
\label{sec:fill}
Although mesh generation is not needed in RBF-FD, certain care still has to be
taken to avoid too close node placements, which cause ill-conditioning. In
other words, we seek a distribution of nodes that conforms to the treated
domain shape and target nodal density, with the nodes also being distributed as
regularly as possible.

\begin{figure}[h!]
  \centering
  \begin{subfigure}[t]{0.4\textwidth}
    \centering
    \hspace{3mm} \includegraphics[width=0.925\linewidth]{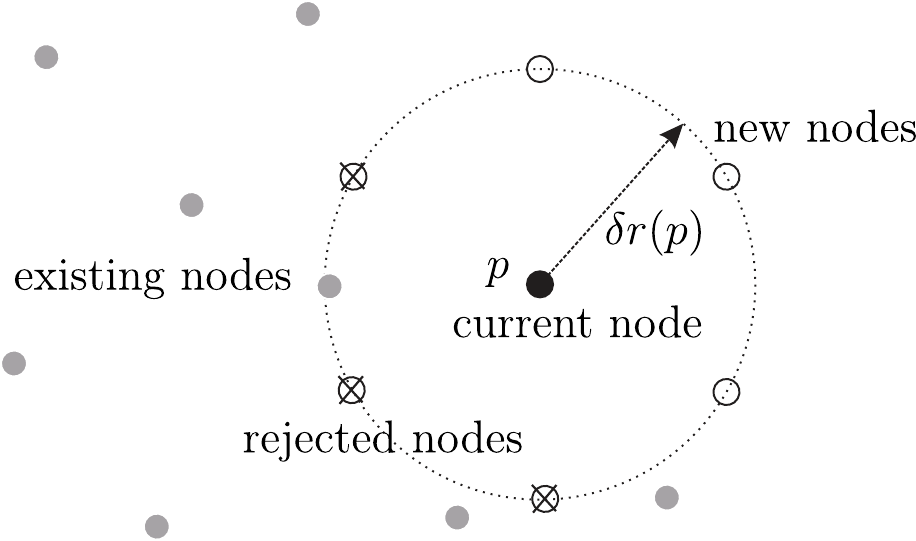}
    \caption{Illustration of the positioning algorithm.}
    \label{fig:poissondisksampling}
  \end{subfigure}
  \begin{subfigure}[t]{0.4\textwidth}
    \centering
    \includegraphics[width=0.925\linewidth]{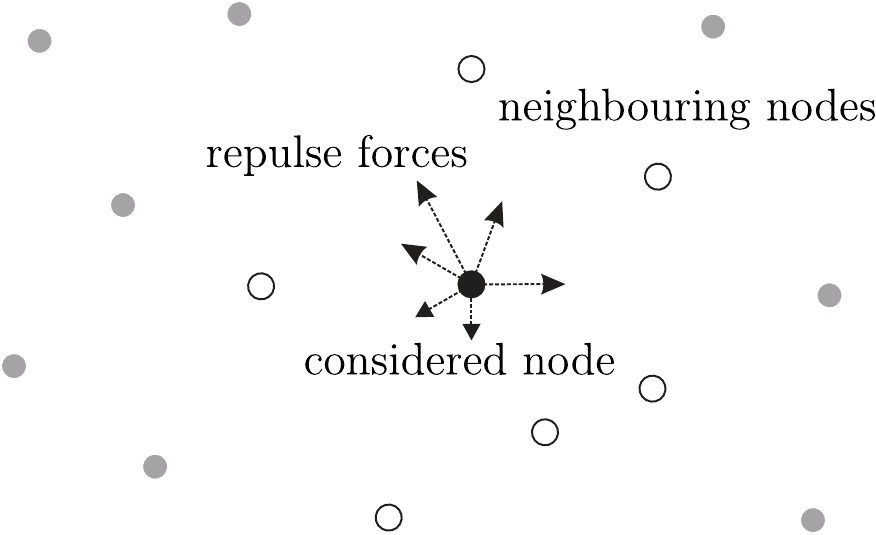}
    \caption{Illustration of the repel algorithm.}
    \label{fig:relax}
  \end{subfigure}
  \caption{Scheme of algorithms used for domain discretisation.}
  \label{fig:positioning}
\end{figure}

The basis for the node placement procedure is the Poisson Disc Sampling
algorithm~\cite{wei2008parallel}. The
procedure begins with adding a starting ``seed node'' into the domain, followed
by an iterative spreading of the populated area of the domain. In each iteration
new nodes are added regularly on a sphere (or circle) with centre in a randomly
selected node $p$ that has not been processed yet and radius $\dr(p)$. The
function $\dr$ is an external parameter and its value at point $p$ represents
the expected distance between neighbouring nodes around $p$, characterising the
nodal density at that point. Some of these newly added nodes might be too close
to the already positioned nodes, and thus a simple rejection criterion is used:
if the distance to the closest node is closer than $\zeta \cdot \dr(p)$, for
$\zeta \in [0, 1]$, the node is removed.  The proximity factor $\zeta$
must be used due to numerical errors in distance computations, however its
use also produces smoother distributions. Is is typically set to $0.9$ or
$0.99$. The algorithm is schematically presented in
Fig.~\ref{fig:poissondisksampling}.

After the domain is populated, the distribution is further regularised
with an iterative algorithm that ``diffuses'' the nodes as if a
repelling force was induced by theirs neighbours. The force
is determined by summing the gradients of potentials in surrounding nodes,
where the amplitude of the force is further normalised by the distance to
the closest neighbour and $\dr$. The algorithm also uses primitive
``simulated annealing'' techniques, i.e.\ the magnitude of translation
is multiplied with an annealing factor that is linearly dropping from $Q_i$ to $Q_f$.
More precisely, in each node a repulsing force is computed as
\begin{equation}
d(i) = -\xi \; \Big( \frac{N_{r} - n_{r} }{N_{r}} \big(Q_i - Q_f) + Q_i\Big)
\sum_{p_j \in \Ns(p_i)} \left(\frac{\dr(p_j)}{\|p_j - p_i\|}\right)^\kappa (p_j - p),
\label{eq:relax}
\end{equation}
where $\xi$ stands for normalization factor, usually distance to the closest
node, $N_r$ represents the number of regularization iterations, $n_r$ current iteration,
$Q_i$ initial ``heat'', $Q_f$ final ``heat'', $\kappa$ potential power and
$\Ns(p)$ are neighbours of $p$. Typical set-up is $N_r \in [3, 10]$, $Q_i \in
[0.8, 1.2]$, $Q_f \in [0, 0.2]$, $\kappa = 2$ and 3 closest neighbours are
considered.  The algorithm is schematically presented in Fig.~\ref{fig:relax}.

\begin{figure}[h!]
  \centering
  \begin{subfigure}[t]{0.31\textwidth}
    \centering
    \includegraphics[width=0.97\linewidth]{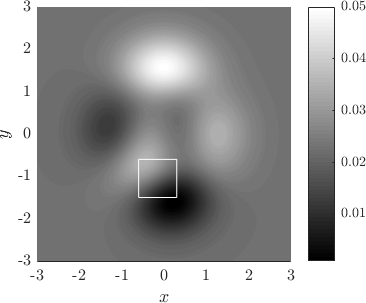}
    \caption{Function $\dr$.}
    \label{fig:sample}
    \vspace{3ex}
  \end{subfigure}
  \begin{subfigure}[t]{0.31\textwidth}
    \centering
    \includegraphics[width=0.78\linewidth]{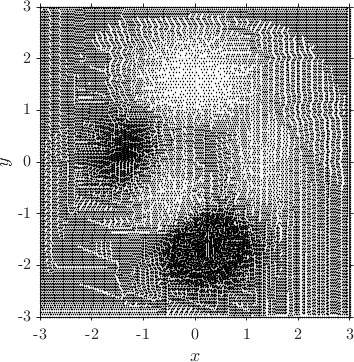}
    \caption{Nodes placed by the fill algorithm.}
    \label{fig:sample-fill}
    \vspace{3ex}
  \end{subfigure}
  \begin{subfigure}[t]{0.31\textwidth}
    \centering
    \includegraphics[width=0.78\linewidth]{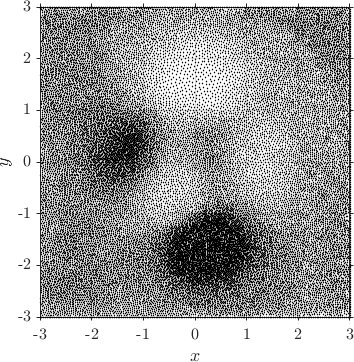}
    \caption{Nodes after regularization.}
    \label{fig:sample-relax}
  \end{subfigure}
  \caption{Illustration of domain discretisation with given function.}
  \label{fig:sample-all}
\end{figure}

\begin{figure}[h!]
  \centering
  \begin{subfigure}[t]{0.31\textwidth}
    \centering
    \includegraphics[width=0.97\linewidth]{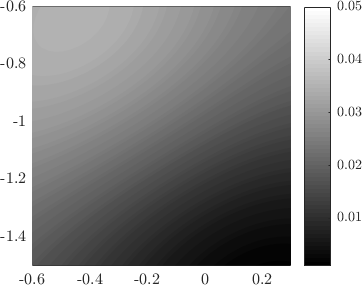}
    \caption{Enlargement of Fig.~\ref{fig:sample}.}
  \end{subfigure}
  \begin{subfigure}[t]{0.31\textwidth}
    \centering
    \includegraphics[width=0.78\linewidth]{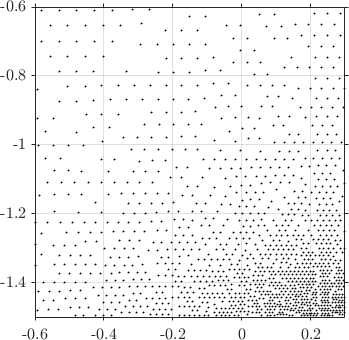}
    \caption{Enlargement of Fig.~\ref{fig:sample-fill}.}
    \label{fig:sample-fill-zoom}
  \end{subfigure}
  \begin{subfigure}[t]{0.31\textwidth}
    \centering
    \includegraphics[width=0.78\linewidth]{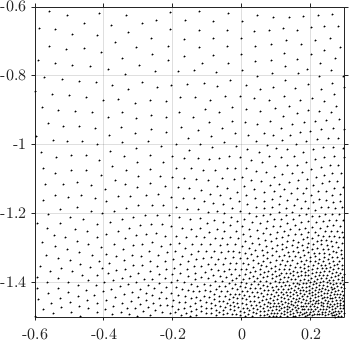}
    \caption{Enlargement of Fig.~\ref{fig:sample-relax}.}
  \end{subfigure}
  \caption{Enlarged versions of plots in Fig.~\ref{fig:sample-all}.}
  \label{fig:sample-all-zoomin}
\end{figure}

A simple demonstration of above algorithms is presented in Fig.~\ref{fig:sample-all}.
A scaled and shifted \texttt{peaks} function from Matlab was used as $\delta r$,
specifically,
\begin{align}
\delta r(x, y) &= \delta x +  (\Delta x - \delta x)(\texttt{peaks}(x, y) + 6.55) / 14.66, \\
\texttt{peaks}(x, y) &= 3(1-x)^2 e^{-x^2 - (y+1)^2} - 10 (x/5-x^2-y^5) e^{-x^2-y^2} - \textstyle\frac13 e^{-(x+1)^2 - y^2}
\end{align}
with $\delta x = 0.007$ and $\Delta x = 0.07$. A square region $[-3, 3] \times [-3, 3]$ was discretised,
filled and regularised with described algorithms. The function $\delta r$ is shown in Fig.~\ref{fig:sample},
and the initial and improved distributions are shown in Figures~\ref{fig:sample-fill}
and~\ref{fig:sample-relax}, respectively. Zoomed-in versions are
presented in Fig.~\ref{fig:sample-all-zoomin} for easier visual assessment of
the node quality.

\begin{figure}[h]
 \centering
 \begin{subfigure}[t]{0.31\textwidth}
   \centering
   \includegraphics[width=0.85\linewidth]{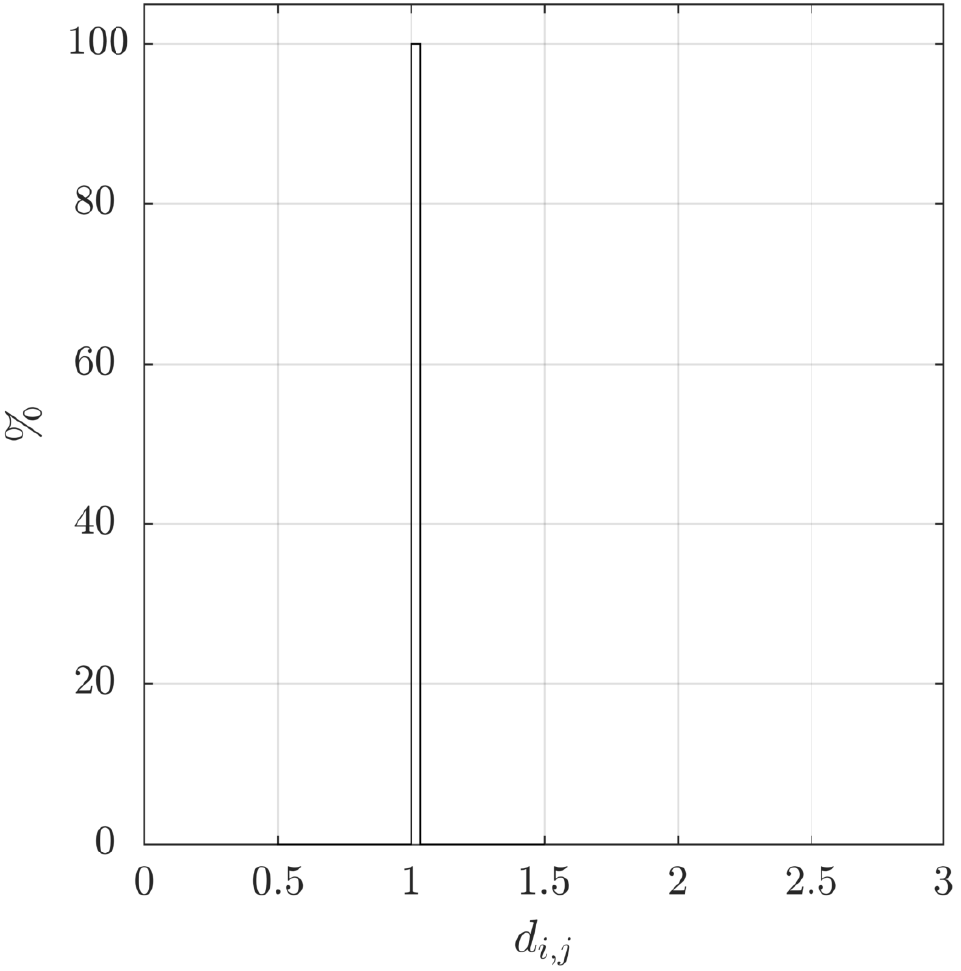}
   \caption{The ``ideal'' histogram.}
 \end{subfigure}
\begin{subfigure}[t]{0.31\textwidth}
  \centering
  \includegraphics[width=0.85\linewidth]{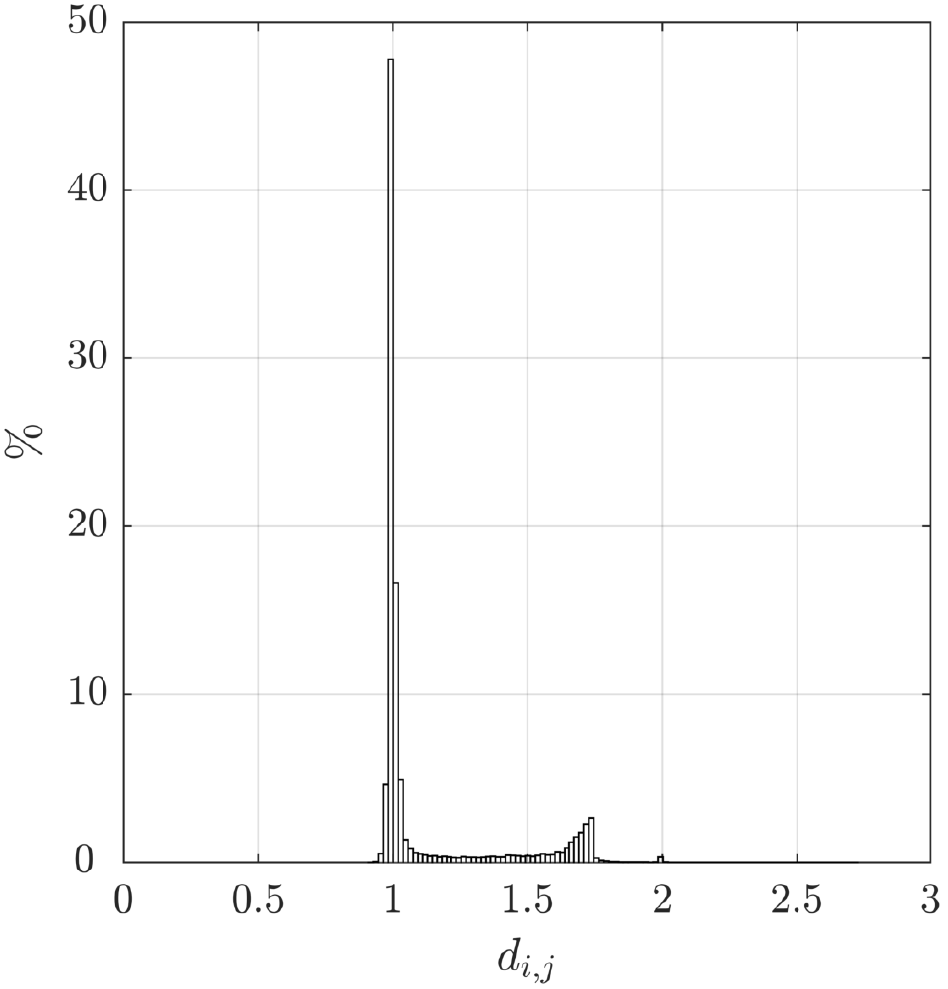}
  \caption{Nodes placed by the fill algorithm.}
\end{subfigure}
\begin{subfigure}[t]{0.31\textwidth}
  \centering
  \includegraphics[width=0.85\linewidth]{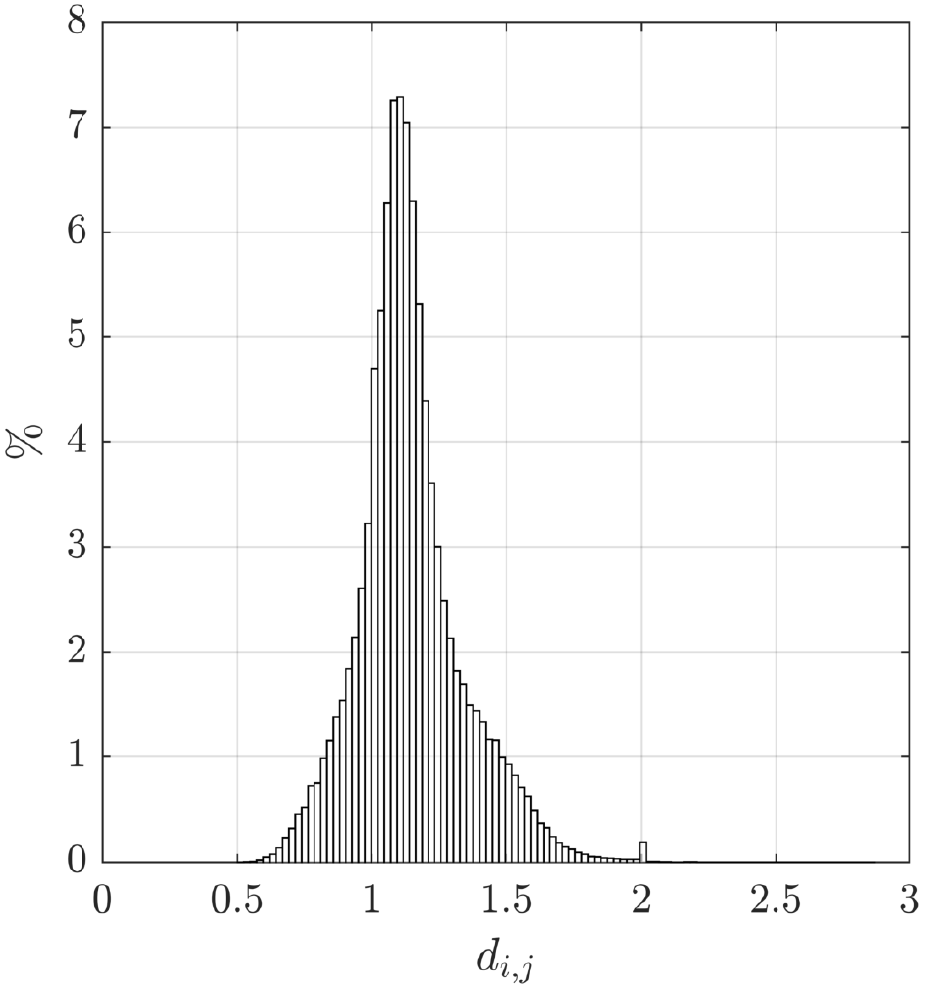}
  \caption{Nodes after regularization.}
\end{subfigure}
\caption{Histogram of normalized distances to 6 closest
nodes.}
\label{fig:dist-hist}
\end{figure}

To present a quantitative measure of the node quality, the distances to the
first layer of closest nodes were computed for each node. Specifically,
if $\{p_{i,j}\}_{j=1}^6$ are the six closest neighbours of $p_i$ (excluding
$p_i$), the normalized distances $d_{i,j}$ were computed as
\begin{equation}
  d_{i,j} = \|p_{i, j} - p_i\| / \dr(p_i).
\end{equation}
Figure~\ref{fig:dist-hist} shows the histograms of normalized distances
$d_{i,j}$ for node distributions before and after regularisation.
The nodal distribution before regularisation is close to the ``ideal''
continuous case, where all $d_{i,j}$ would be equal to 1. The tail
peak of this histogram at around $1.7$ has been investigated and
corresponds to the nodes at the gaps in the distribution.
A small portion of $d_{i,j}$ at $d_{i,j} = 2$ correspond to the nodes
at the boundary where some of the six nodes fall in the second
layer of closest nodes. The histogram of the regularised distribution
is smoother with a peak at $d_{i,j} = 1.2$ still corresponding to the nodes at
the boundary. The peak at $1.7$ before the regularisation, which corresponded to
the gaps in the distribution, has been removed at the cost of greater spread.

Node generation is relatively undemanding from the computational point of view. The
most expensive part of the fill algorithm is determining whether there are any proximity
violations, which requires searching the closest nodes. However, this can be optimized
with appropriate dynamic tree structures and managing the candidate list of
possible nearby nodes. The regularization algorithm is a bit more expensive since it
requires a full neighbourhood search in each iteration, but typically few iterations
are needed to smooth out the irregularities.

Further analyses of the filling algorithm are out of the scope of this paper, but
an interested reader is referred to a paper~\cite{slak2018generation},
which provides a more quantitative analysis of the fill algorithm along with comparison
to other existing techniques.

\subsection{Local density modification}
\label{sec:adaptivity}
\begin{figure}[h]
  \centering
  \includegraphics[height=5cm]{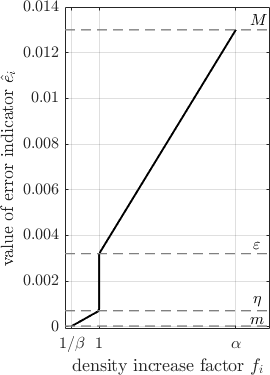}
  \includegraphics[height=5cm]{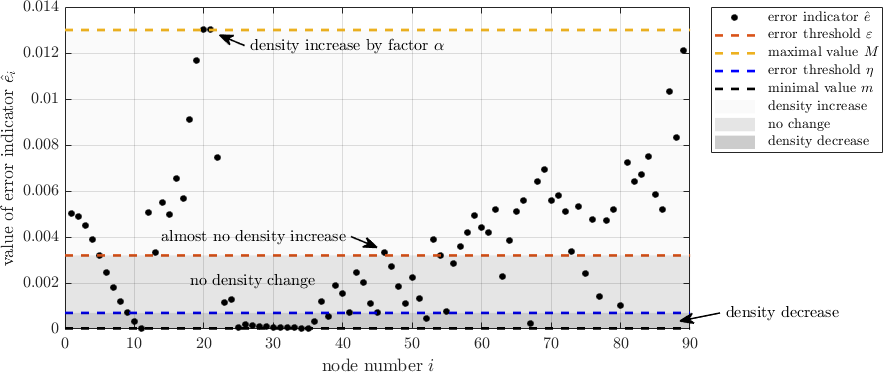}
  \caption{Construction of new density function $\widetilde{\dr}$ from $\dr$ and $\eh$.}
  \label{fig:new-density}
\end{figure}

After the numerical solution of the PDE has been computed and the error indicator computed,
the local nodal density is altered.
This is achieved by computing new nodal radii $\widetilde{\dr}_i$ from the existing ones as
\begin{equation}
\widetilde{\dr}_i = \min(\dr_i / f_i, \drup_i),
\label{eq:adapt-rho}
\end{equation}
where $\dr_i$ is the distance to from $p_i$ to the closest neighbouring node.
This imposes the upper bound on radii from point~(\ref{enum:deref-bnd}) from the preliminaries above.
Similarly, a lower bound $\drlo$ limiting refinement could also be imposed,
however in practice that bound was not needed. The factor $f_i$ is called a \emph{density increase factor}
and represents the relative change in nodal density at point $p_i$.
It is computed as
\begin{equation}
\label{eq:adapt-rho-factor}
f_i = \begin{cases}
1 + \frac{\eta - \eh_i}{\eta - m} (\frac{1}{\beta} - 1), & \eh_i \leq \eta, \quad \text{i.e.\ decrease the density} \\
1, & \eta < \eh_i < \eps,  \quad \text{i.e.\ no change in density}\\
1 + \frac{\eh_i - \eps}{M - \eps} (\alpha - 1), & \eh_i \geq \eps, \quad \text{i.e.\ increase the density}
\end{cases}
\end{equation}
where $M = \max_{p_i \in \P} \eh_i$ is the maximal and $m = \min_{p_i \in \P} \eh_i$
the minimal value of the error indicator. The plot of factor $f_i$ with respect to
the error indicator $\eh_i$ can be seen in Fig.~\ref{fig:new-density} (left).

Note that $\frac{1}{\beta} \leq f_i \leq \alpha$ always holds, as
$\frac{\eh_i - \eps}{M - \eps} \in [0, 1]$,
and $\alpha \geq 1$. If the node is on the error threshold, i.e.~$\eh_i = \eps$,
the factor $f_i$ equals $1$ and the density will stay the same,
ensuring compatibility with the case when $\eh_i < \eps$.
Additionally, if the node has the highest error, i.e.~$\eh_i = M$,
the increase factor will be maximal, i.e.\ $\alpha$.
Symmetric observations hold for the derefine case. Setting
$\alpha = 1$ or $\beta = 1$ disables refinement or derefinement, respectively.

Once the desired nodal radii $\widetilde{\dr}_i$ have been evaluated in all nodes,
a new continuous function  $\widetilde{\dr}$ is reconstructed using
Modified Sheppard's method as described in the next section.
The overall density modification procedure is illustrated conceptually in Fig.~\ref{fig:new-density}.

\subsection{Computing the new density function}
\label{sec:reconstruct-density}
After new desired local radii $\widetilde{\dr}_i$~\eqref{eq:adapt-rho} have been computed,
a continuous density function suitable for the fill algorithm in section~\ref{sec:fill}
needs to be constructed. To construct a function $\dr(p)$ from discrete values
$\widetilde{\dr}_i$, the Modified Sheppard's method is used.  In a special case of
$\alpha=\beta=1$ this becomes an inverse problem of the node positioning problem described in
section~\ref{sec:fill}.
The problem of computing the desired internodal distance at an arbitrary point
$p \in \P$ translates to a problem of scattered data interpolation with known
values $\widetilde{\dr}_i$ at points $p_i$. Standard methods for this are well developed,
described and compared by Franke~\cite{franke1982scattered}. Modified Shepard's
method~\cite[p.\ 185]{franke1982scattered} was used in this paper, where value
of $\dr$ at a point $p$ is computed as a weighted sum of $n$ closest neighbours
of $p$
\begin{align}
\dr(p) = \begin{cases}
\widetilde{\dr}_i, & \text{if } p = p_i \text{ for some } i \\
\frac{\sum_{i=1}^{n} w_i(p)\widetilde{\dr}_i} {\sum_{i=1}^n w_i(p)} & \text{else}
\end{cases}
\end{align}
and the inverse distance weights $w_i(p)$ are defined by
\begin{equation}
w_i(p) = \left((1 - \textstyle \frac{d_i}{d_n}) \frac{1}{d_i} \right)^2,
\end{equation}
where $d_i$ represents the Euclidean distance from $p$ to its $i$-th
neighbour and the neighbours are sorted according to $d_i$, thus making
$d_n$ the largest.

There is however no need for interpolation of values $\widetilde{\dr}_i$, as they themselves
are only an approximation, and techniques from scattered data
approximation~\cite{wendland2004scattered} might be more appropriate. We have
tested different approximation and interpolation methods for values $\dr_i$ with
no significant changes in characteristics of node distributions upon refilling
the domain.  Modified Sheppard's method was chosen based on computational
efficiency and ease of implementation.

\subsubsection{Analysis of repeated fill and reconstruct cycles}

\begin{figure}[h!]
  \begin{subfigure}[t]{0.31\textwidth}
    \centering
    \includegraphics[width=4cm]{images/sample_relaxed_full.png}
    \caption{Initial nodes.}
    \label{fig:sample-fill-2}
    \vspace{3ex}
  \end{subfigure}
  \begin{subfigure}[t]{0.33\textwidth}
    \centering
    \includegraphics[width=4cm]{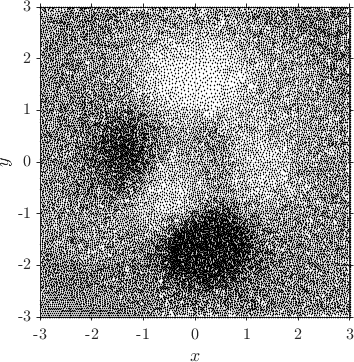}
    \caption{Nodes after 1 fill-reconstruct cycle.}
    \label{fig:sample-fill-2d-iter}
    \vspace{3ex}
  \end{subfigure}
  \begin{subfigure}[t]{0.35\textwidth}
    \centering
    \includegraphics[width=4cm]{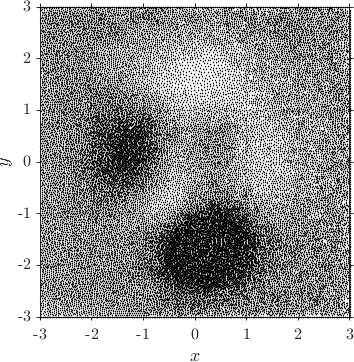}
    \caption{Nodes after 10 fill-reconstruct cycles.}
    \label{fig:sample-repeat}
  \end{subfigure}
  \caption{Analysis of fill-reconstruct cycles.}
  \label{fig:reconstruct}
\end{figure}

Every iteration of the adaptive procedure consists of filling the domain,
solving the problem, and reconstructing and adapting the node density.
Even without adaptation, these fill-reconstruction
cycles do not guarantee the preservation of nodal density and it is therefore
important to analyse their behaviour before implementing this approach into the
refinement algorithm. To that end a continuation of the test case from
section~\ref{sec:fill} was performed to better asses the effects of the repeated application
of fill-reconstruction cycles. The effects of 10 successive cycles
on the initial nodal distribution are demonstrated in Fig.~\ref{fig:sample-repeat}.

Repeated application of fill-reconstruction cycles tends to slightly increase the nodal density, due
to the proximity factor that is set to $\zeta=0.99$.  After 1 iteration, the
minimal internodal distance is reduced from $0.00699727$ to $0.00697768$ and after
10 iterations to $0.00688613$, i.e.\ a reduction of 1.5\% over the course of 10
iterations, which is well within the allowed limit of 1\% per iteration, as
imposed by the fill tolerance factor $\zeta$. The most visible effect of
successive application is the gradual spread and smear of high density areas,
which expand with every successive iteration. This feature would be undesirable
in certain situations and can be dealt with, however in numerical examples
presented in this paper the phenomenon did not have considerable impact and
is therefore not discussed further.

\subsection{Algorithm}

The presented adaptive procedure is summarized as Algorithm~\ref{alg:adapt}.
The \textsc{discretize} procedure refers to the fill algorithm described in
section~\ref{sec:fill}, the \textsc{solve} procedure refers to the application of the RBF-FD
method for numerical solution of PDEs as described in section~\ref{sec:rbffd} and the
\textsc{adapt} procedure refers to the local density modification and function reconstruction
as described in sections~\ref{sec:adaptivity} and~\ref{sec:reconstruct-density}.

\begin{algorithm}[h]
  \caption{Adaptive solution procedure}
  \label{alg:adapt}
  \vspace{1mm}
  \textbf{Input:} The problem, computation domain $\Omega$, initial density function $\dr\colon\Omega\to\R$, the maximal number of iterations $I_{\text{max}}$ and adaptivity parameters $\eps, \eta, \alpha, \beta, \drup$. \\
  \textbf{Output:} The numerical solution of the problem.
  \begin{algorithmic}[1]
    \Function{adaptive\_solve}{problem, $\Omega, \dr, I_{\text{max}}, \eps, \eta, \alpha, \beta, \drup$}
    \For{$i$}{0}{$I_{\text{max}}$} \label{alg:while}
    \State $\P \gets \Call{discretise}{\Omega, \dr}$
    \Comment{Discretises domain $\Omega$ as described in section~\ref{sec:fill}.}
    \State $\text{solution} \gets \Call{solve}{\text{problem}, \P}$
    \Comment{Solves the problem using discretisation $\P$, see section~\ref{sec:rbffd}.}
    \State $\text{error} \gets \Call{estimate\_error}{\text{solution}, \P}$
    \Comment{Error indicator computation, see~\eqref{eq:anal-err-est} or~\eqref{eq:err-indicator-std}.}
    \If{$\Call{mean}{\text{error}} \geq \eps$}
    \Comment{Other error reductions such as \textsc{Max} can be used.}
    \State \Return $\text{solution}$
    \EndIf
    \State $\dr \gets \Call{adapt}{\dr, \P, \eps, \eta, \alpha, \beta, \drup}$
    \Comment{Adapt the discretisation as described in section~\ref{sec:adaptivity}.}
    \EndFor
    \State \Return $\text{solution}$
    \EndFunction
  \end{algorithmic}
\end{algorithm}

\subsection{An illustrative adaptivity example}
To help illustrate the proposed adaptive framework we use it on a simple test case of
function approximation. The function
\begin{equation}
  g(x) = 3(1-x)^2 \exp(-x^2) + 3\exp(-4(x-1)^2)
\end{equation}
is approximated on a given set of points using moving weighted least squares with
basis $\{1, x\}$, 12 support nodes and Gaussian weight $\exp(-x^2/d^2)$, where $d$ is the distance to
the farthest support node. The resulting approximation is denoted $\hat{g}$.

\begin{figure}[h]
  \centering
  \includegraphics[width=0.4\linewidth]{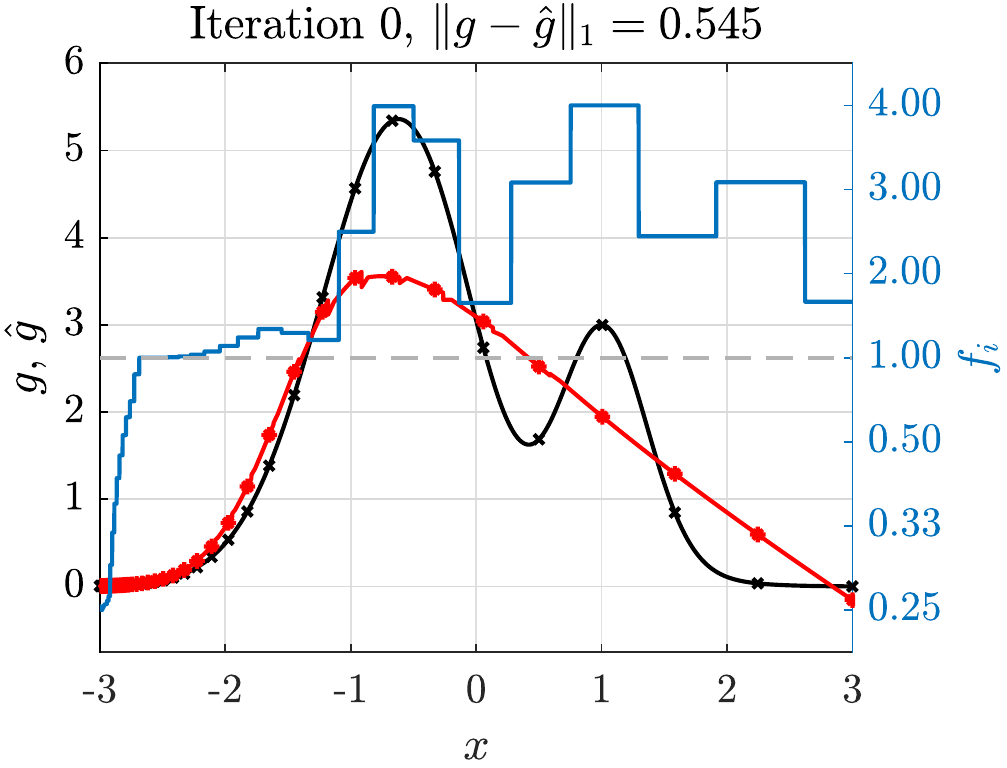}
  \includegraphics[width=0.4\linewidth]{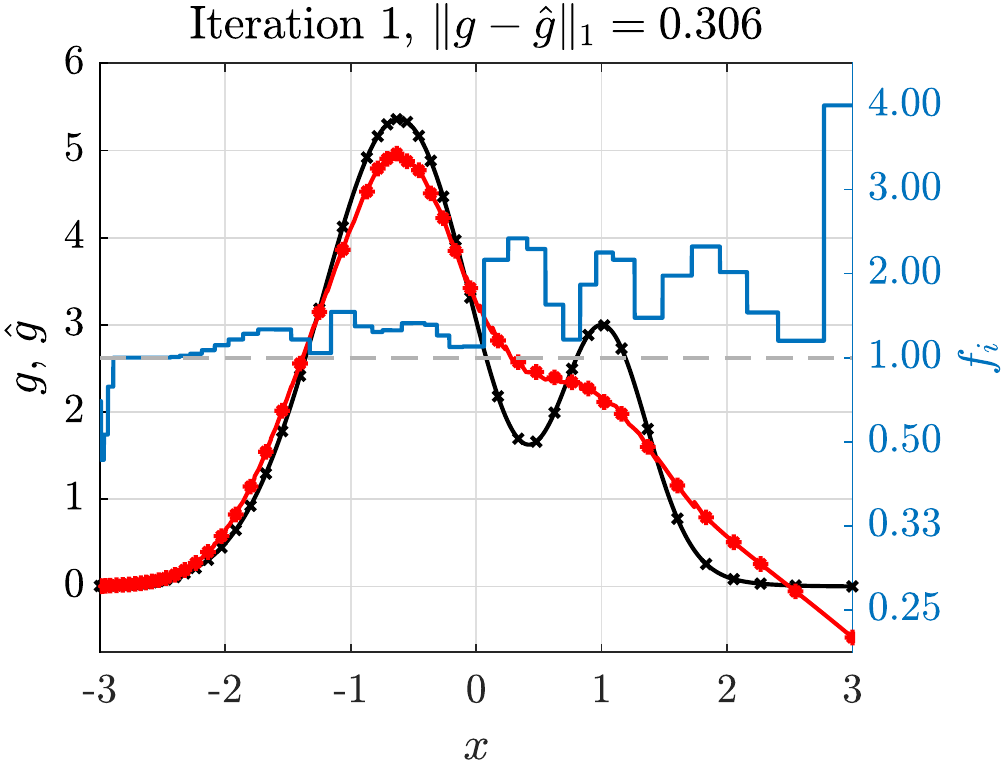} \\
  \includegraphics[width=0.4\linewidth]{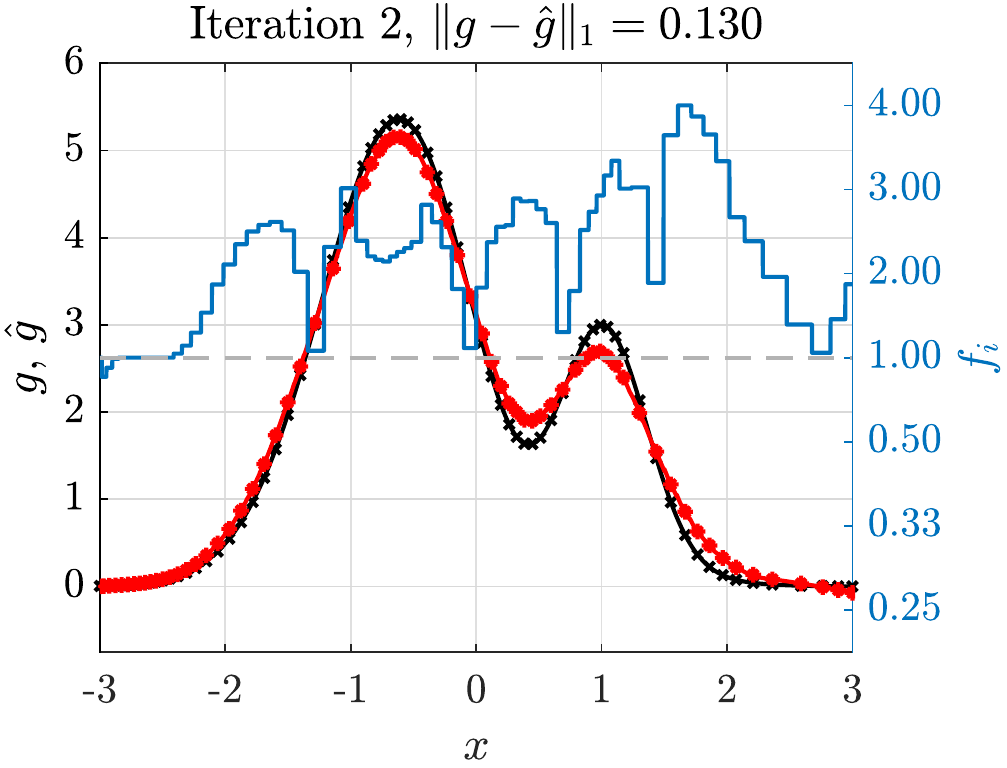}
  \includegraphics[width=0.4\linewidth]{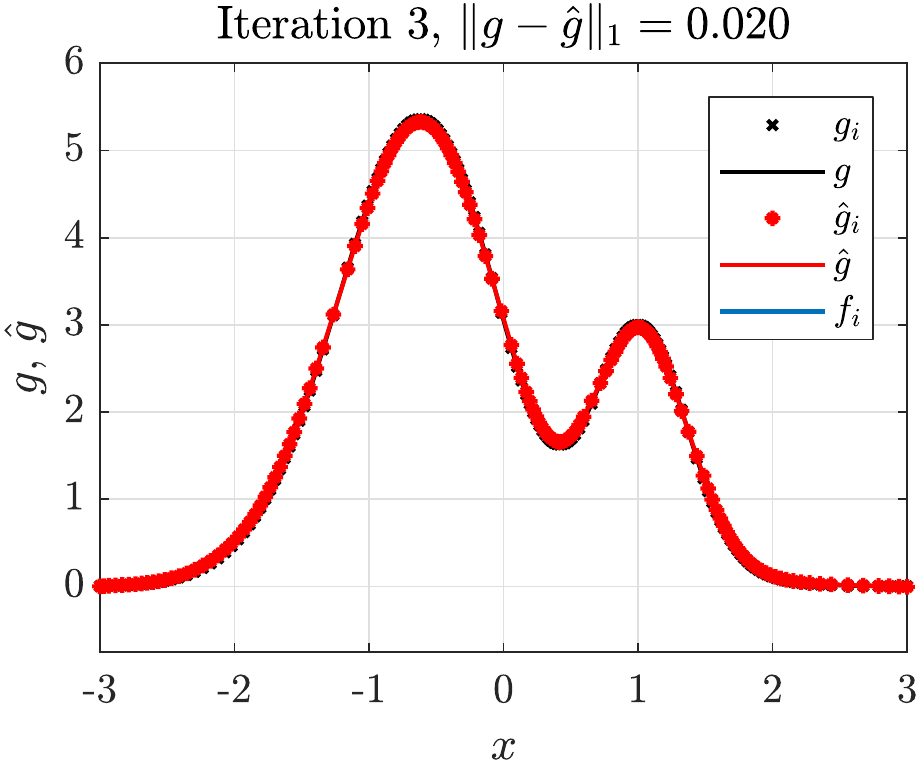}
  \caption{Illustration of adaptive function approximation using Algorithm~\ref{alg:adapt}.}
  \label{fig:adaptive-illustration}
\end{figure}

Initially, the domain is filled with density
\begin{align}
  \dr(x) = h_0(1+25|3+x|), \quad h_0 = 0.005,
\end{align}
which results in a dense distribution on the negative part of the $x$-axis and a coarse distribution on
the positive part (see Fig.~\ref{fig:adaptive-illustration}, iteration 0). The upper bound for $\dr$ is set
to $\drup(x) = 10 h_0$. The error is computed simply by taking the absolute difference between the
actual and approximation value, multiplied by the local nodal spacing,
$e_i = |g(x_i) - \hat{g}(x_i)| \dr(x_i)$.
Fig.~\ref{fig:adaptive-illustration} shows 3 iterations of the Algorithm~\ref{alg:adapt}
with error thresholds $\eps = 10^{-3}$, $\eta = 10^{-4}$, and aggressivenesses $\alpha = 4$ and $\beta = 4$.

Initially, the function $\hat{g}$ approximates $g$ well around $x=-3$ and poorly for $x > -1$, which is also
reflected in the density increase factor $f_i$. It is greater than 1 and proportional to the error on the
positive side, reaching a maximum of $\alpha$, where the error is greatest. Around $x=-3$ the function
$g$ is approximated well, the error is below $\eta$ and $f_i$ is less than one and dropping to $1/\beta$,
indicating derefinement is required. The node distribution in iteration $1$ is adapted accordingly,
with denser areas around peaks. This iteration continues and the $L_1$ error $\|g-\hat{g}\|_1$,
shown in figure titles, decreases accordingly. The number of nodes, along with the counts of refined and
derefined nodes, is shown in Table~\ref{tab:illustr-n}.

\begin{table}[h]
  \centering
  \caption{Number of nodes during the course of adaptive iteration shown in
  Fig.~\ref{fig:adaptive-illustration}.}
  \label{tab:illustr-n}
  \begin{tabular}{r|r|r|r|r|l}
    iteration & total nodes & refined & no change & derefined & derefined, but hit bound $\drup$ \\   \hline
    \rule{0pt}{12pt} 0 & 41 & 18 & 4 & 11 & 8 \\
    1 & 50 & 38 & 8 & 2 & 2 \\
    2 & 72 & 59 & 10 & 1 & 2 \\
    3 & 159 & / & / & / & / \\
  \end{tabular}
\end{table}

\section{Numerical results}
\label{sec:results}

\subsection{Governing problem}
The solution procedure presented in previous section is applied to linear
elasticity contact problems, more precisely to problems of deformation of elastic
homogeneous isotropic materials under load. The problem is governed by
Cauchy-Navier equations
\begin{equation}
(\lambda + \mu) \nabla (\nabla \cdot \vec{u}) + \mu \nabla^2 \vec{u} = \vec{f},
\end{equation}
where $\vec{u}$ are unknown displacements, $\vec{f}$ is the loading body force,
and $\lambda$ and $\mu$ are Lam\'{e} parameters, often expressed in terms of
Young's modulus $E$ and Poisson's ratio $\nu$.  Another important quantity is
the stress tensor $\sigma$, related to strain via Hooke's law
\begin{equation}
\label{eq:elasticity}
\sigma = \lambda \operatorname{tr}(\eps) I + 2 \mu \eps, \quad \eps = \frac{\nabla \vec{u} + (\nabla \vec{u})^\T}{2},
\end{equation}
where $\lambda$ and $\mu$ are Lam\'{e} parameters from above and $I$ is the
identity tensor.

Two types of boundary conditions are considered: essential boundary conditions
which specify displacements on some portion of the boundary of the domain, i.e.\
$\vec{u} = \vec{u}_0$, and traction boundary conditions, which prescribe surface
traction $\sigma\vec{n} = \vec{t}_0$, with $\vec{n}$ being an outside unit
normal to the boundary of the domain.

In two dimensions, we will use the following componentwise notation
for $\vec{u}$ and $\sigma$, for sake of simplicity:
\begin{equation}
\vec{u} = (u, v) \quad \text{and} \quad
\sigma = \begin{bmatrix} \sigma_{xx} & \sigma_{xy} \\
\sigma_{xy} & \sigma_{yy} \end{bmatrix}.
\end{equation}

\subsection{Disk under pressure}
\label{sec:half2d}
\begin{figure}[h]
  \centering
  \includegraphics[width=0.25\linewidth]{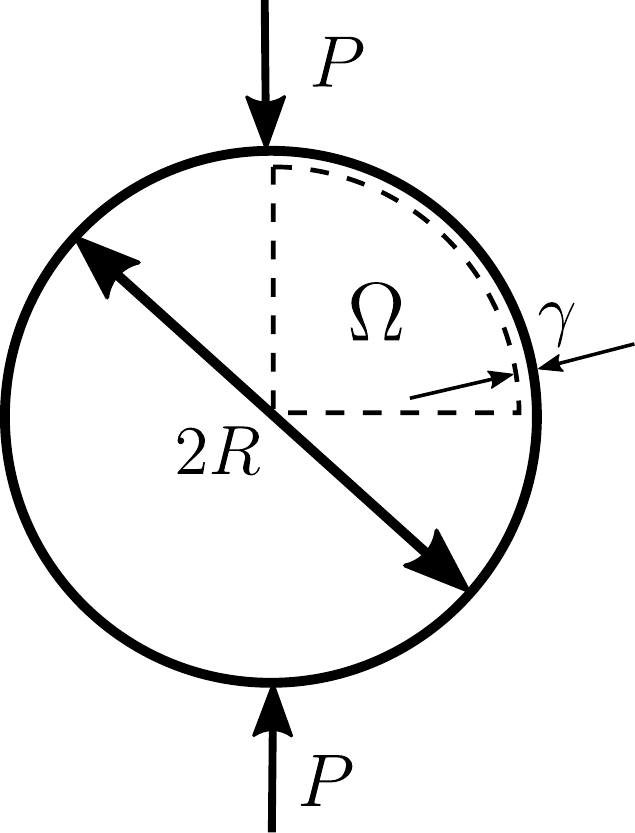}
  \caption{Disk under diametrical compression.}
  \label{fig:disk}
\end{figure}

To demonstrate the adaptive strategy from section~\ref{sec:adaptivity},
a 2D example from linear elasticity, inspired by~\cite{afshar2011node},
is chosen. A disk of radius $R$ is subjected to diametrical compression
with a point force of magnitude~$P$, as illustrated in Fig.~\ref{fig:disk}.

Exact closed form solution for stresses, assuming plane stress condition, is
known, given by e.g.\ Sadd~\cite[p.\ 197]{sadd2009elasticity} as:
\begin{align}
\label{eq:half-disk-anal-sxx}
\sigma_{xx} &= -\frac{2 P}{\pi } \left(\frac{x^2 (R-y)}{r_1^4}+\frac{x^2 (R+y)}{r_2^4}-\frac{1}{2 R}\right) \\
\label{eq:half-disk-anal-syy}
\sigma_{yy} &= -\frac{2 P}{\pi} \left(\frac{(R-y)^3}{r_1^4}+\frac{(R+y)^3}{r_2^4}-\frac{1}{2 R}\right) \\
\label{eq:half-disk-anal-sxy}
\sigma_{xy} &= \frac{2 P}{\pi} \left(\frac{x (R-y)^2}{r_1^4}-\frac{x(R+y)^2}{r_2^4}\right)
\end{align}
where $r_{1,2} = \sqrt{x^2 + (R\mp y)^2}$ are the distances to the poles.
The stress distribution is singular at the poles.

Numerically, only a quarter of the disk is considered due to symmetry.
Additionally, a small boundary layer of width $\gamma$ is subtracted from the
curved edge of the domain, to avoid the singularity at the poles. Therefore, the computational domain
$\Omega$, shown in Fig.~\ref{fig:disk}, is a quarter of a circle with radius
$R - \gamma$. The parameter $\gamma$ determines the peak stress magnitudes of
order $\frac{P}{2\pi \gamma}$ concentrated around the point $x=0$, $y=R-\gamma$,
allowing fine control over case difficulty. Thus the following problem
is solved:
\begin{align}
&(\lambda + \mu) \nabla (\nabla \cdot \vec{u}) + \mu \nabla^2 \vec{u} = 0 \quad \text{ in } \quad
\Omega = \{(x, y) \in \R^2;\ x \geq 0, y \geq 0, x^2+y^2 \leq (R-\gamma)^2 \}.
\label{eq:compressed-disk-problem}
\end{align}
Traction boundary conditions are imposed on the curved boundary,
given by equations~(\ref{eq:half-disk-anal-sxx}--\ref{eq:half-disk-anal-sxy}).
Symmetry boundary conditions $v = 0, \dpar{u}{y} = 0$ and $u = 0, \dpar{v}{x} = 0$
are imposed on the bottom and left boundary, respectively.

\begin{figure}[h]
  \centering
  \includegraphics[width=\linewidth]{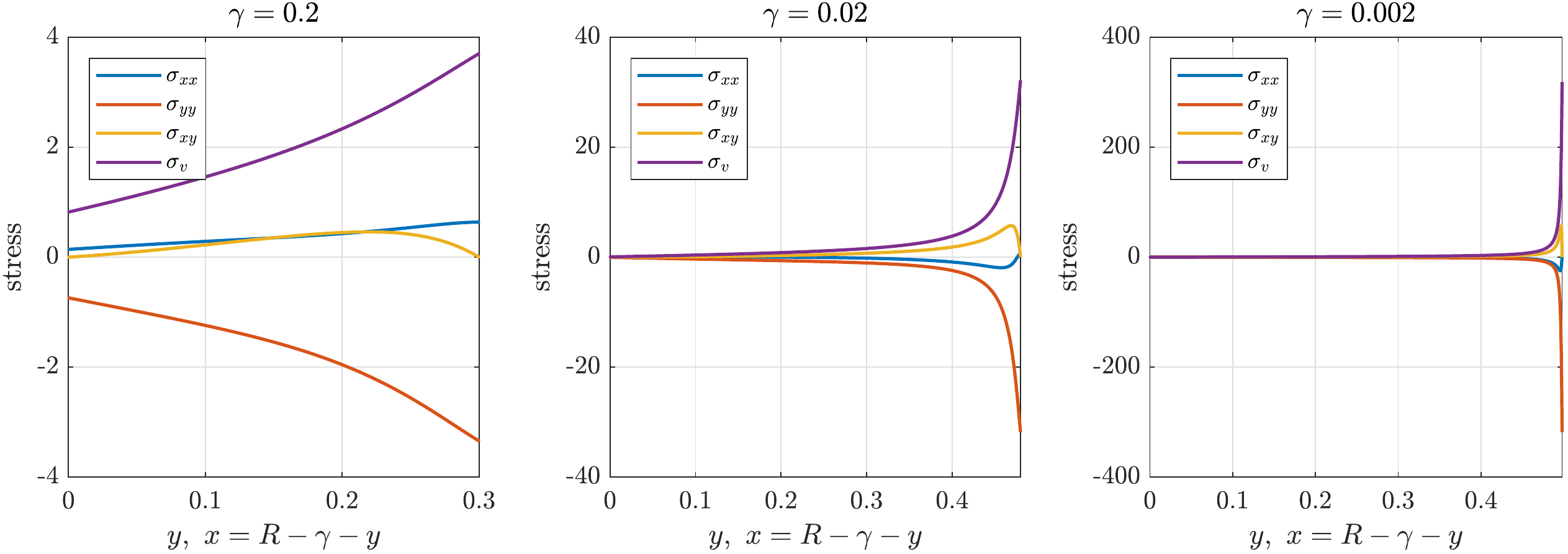}
  \caption{Stress values along the antidiagonal $x = R-\gamma-y$ for three different cases.}
  \label{fig:anal-cross}
\end{figure}

Three possible values for $\gamma$ will be considered, $\gamma_1 = 0.2$,
$\gamma_2 = 0.02$ and $\gamma_3 = 0.002$. The values of the stresses along the
antidiagonal $x = R-\gamma-y$ of $\Omega$ are shown in
Fig.~\ref{fig:anal-cross}.

The error between analytical stresses $\sigma$ and computed stresses $\hat{\sigma}$
is measured in three different ways:
\begin{itemize}
  \item Relative approximate $L^\infty$ norm:
  \begin{equation}
    \label{eq:errinf}
  e_\infty = \frac{\|\hat{\sigma} - \sigma\|_\infty}{\|\sigma\|_\infty}, \quad
  \left\|\sigma \right\|_\infty = \max_{p \in \Omega_h} \max\{ \sigma_{xx}(p), \sigma_{yy}(p), \sigma_{xy}(p)\}
  \end{equation}
  \item Relative approximate $L^1$ norm:
  \begin{equation}
  \label{eq:err1}
  e_1 = \frac{\|\hat{\sigma} - \sigma\|_1}{\|\sigma\|_1}, \quad
  \left\|\sigma \right\|_1 = \frac{1}{3 |\Omega_h|} \sum_{p \in \Omega_h} (|\sigma_{xx}(p)| + |\sigma_{yy}(p)| + |\sigma_{xy}(p)|)
  \end{equation}
  \item Relative energy norm:
  \begin{equation}
  e_E = \frac{\|\hat{\sigma} - \sigma\|_E}{\|\sigma\|_E}, \quad
  \left\|\sigma \right\|_E^2 = \int_\Omega \tilde{e}_E(p) \, d\Omega, \quad \tilde{e}_E(p) = \sigma(p) : S : \sigma (p),
  \label{eq:energy-norm}
  \end{equation}
  where $S$ is the fourth-order compliance tensor computed from $E$ and $\nu$.
\end{itemize}

The norms are evaluated by discretising $\Omega$ with a dense uniform grid
$\Omega_h$ with grid spacing $h = 5\cdot 10^{-4}$, resulting in around $800\,000$
grid points. For $e_\infty$ and $e_1$ the maximum and the sum are taken over
these grid points. The intermediate values of $\sigma(p)$ for $p \notin \P$
were evaluated using linear interpolation over the Delaunay triangulation
induced by the points in $\P$. The integral in $e_E$ was approximated using
trapezoidal rule over $\Omega_h$.

Adaptive approach, described in section~\ref{sec:adaptivity} is now employed to
solve problem~\eqref{eq:compressed-disk-problem}. The error indicator is
constructed from the known analytical
solution~(\ref{eq:half-disk-anal-sxx}--\ref{eq:half-disk-anal-sxy}), using the
kernel of the integral for energy norm~\eqref{eq:energy-norm} multiplied by the
approximation of the volume element
\begin{equation}
\label{eq:anal-err-est}
\hat{e}(p) = (\sigma(p)\!:S\!:\sigma(p)) \dr(p)^d,
\end{equation}
where $d=2$ is the dimension of the problem.

To discretise the domain, fill algorithm from section~\ref{sec:fill} with $\zeta
= 0.9$ was used and the distribution was regularised using repel algorithm with
$N_r = 10$ iterations, $Q_i = 0.8$, $Q_j = 0$, $\kappa=2$ and 3 neighbours.
The domain was initially filled with constant density of $\drI = 0.02$, amounting to
approximately $N = 600$ nodes in the domain.  For the given discretisation, RBF-FD
with 25 neighbouring nodes and 15 Gaussian RBFs was used. The shape parameter
$\sigma_b$ was varied proportionally to the internodal distance with base value
$\sigma_b = 100$.  For adaptive iteration, $\eps = 10^{-7}$, $\eta = 10^{-9}$,
$\alpha = 5$ and $\beta = 1.5$ were chosen. The reconstruction and evaluation of
the new density function $\widetilde{\dr}$ was done using $7$ closest
neighbours.  Values of $R=0.5$, $P=1$, $E=1$ and $\nu = 0.33$ were taken for
physical parameters of the problem.

The overall procedure seems to depend on a lot of free parameters, however most are
either well studied, such as the RBF shape parameters (see
section~\ref{sec:rbffd-shape-comp}) or their value is not very important to the
overall procedure, such as the fill, repel and reconstruct parameters, described
in more detail in sections~\ref{sec:fill} and~\ref{sec:reconstruct-density}.
Standard values, as described above and as used by Slak and Kosec~\cite{slak2018refined}, were
taken for these parameters and they were kept the same for all runs, to avoid
hyper-tuning the results. The parameters related to adaptivity are further
analysed in section~\ref{sec:adaptivity-params} and justify the choices above.

\begin{figure}[h!]
  \centering
  \includegraphics[width=0.25\linewidth]{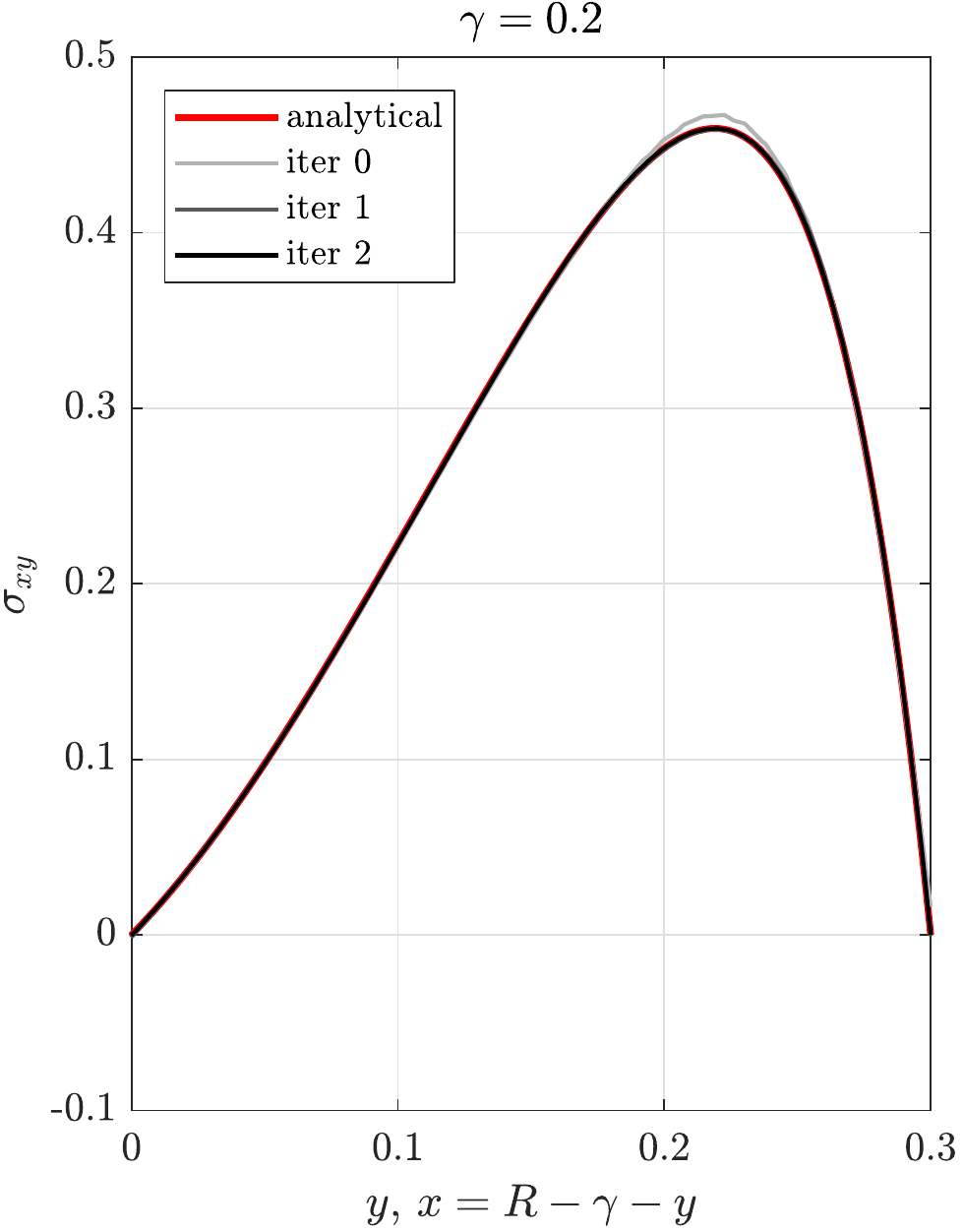}
  \includegraphics[width=0.25\linewidth]{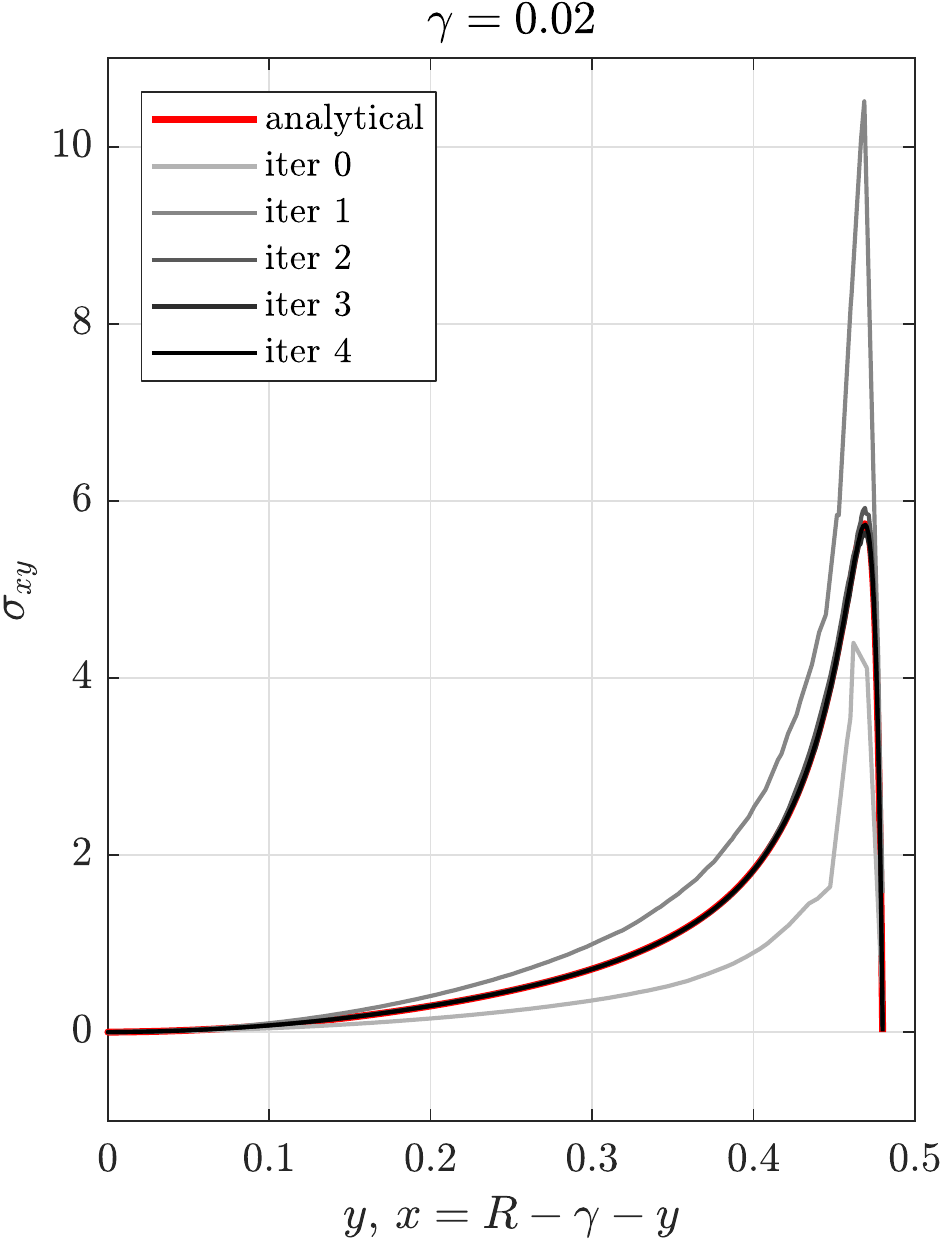}
  \includegraphics[width=0.25\linewidth]{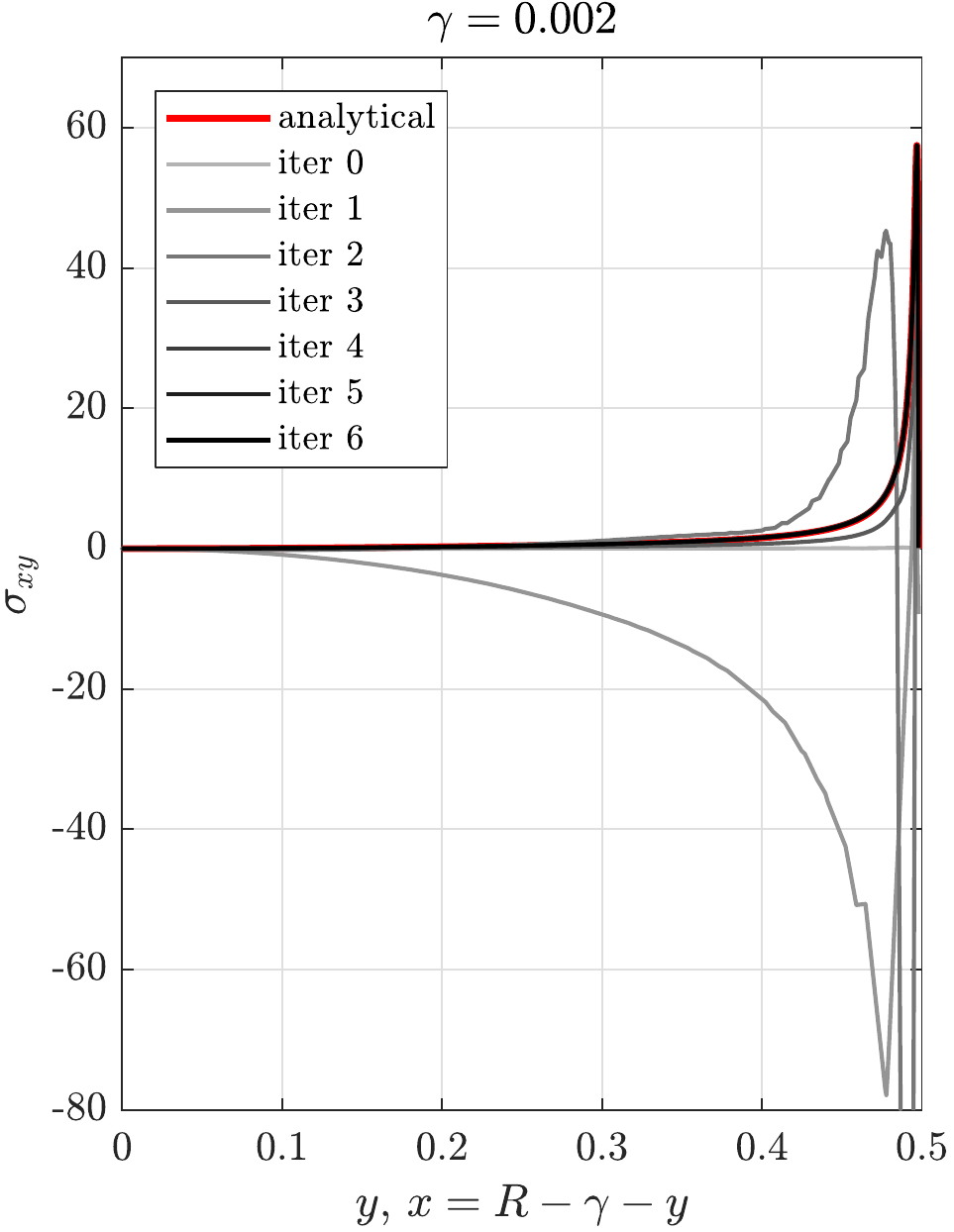}

  \includegraphics[width=0.25\linewidth]{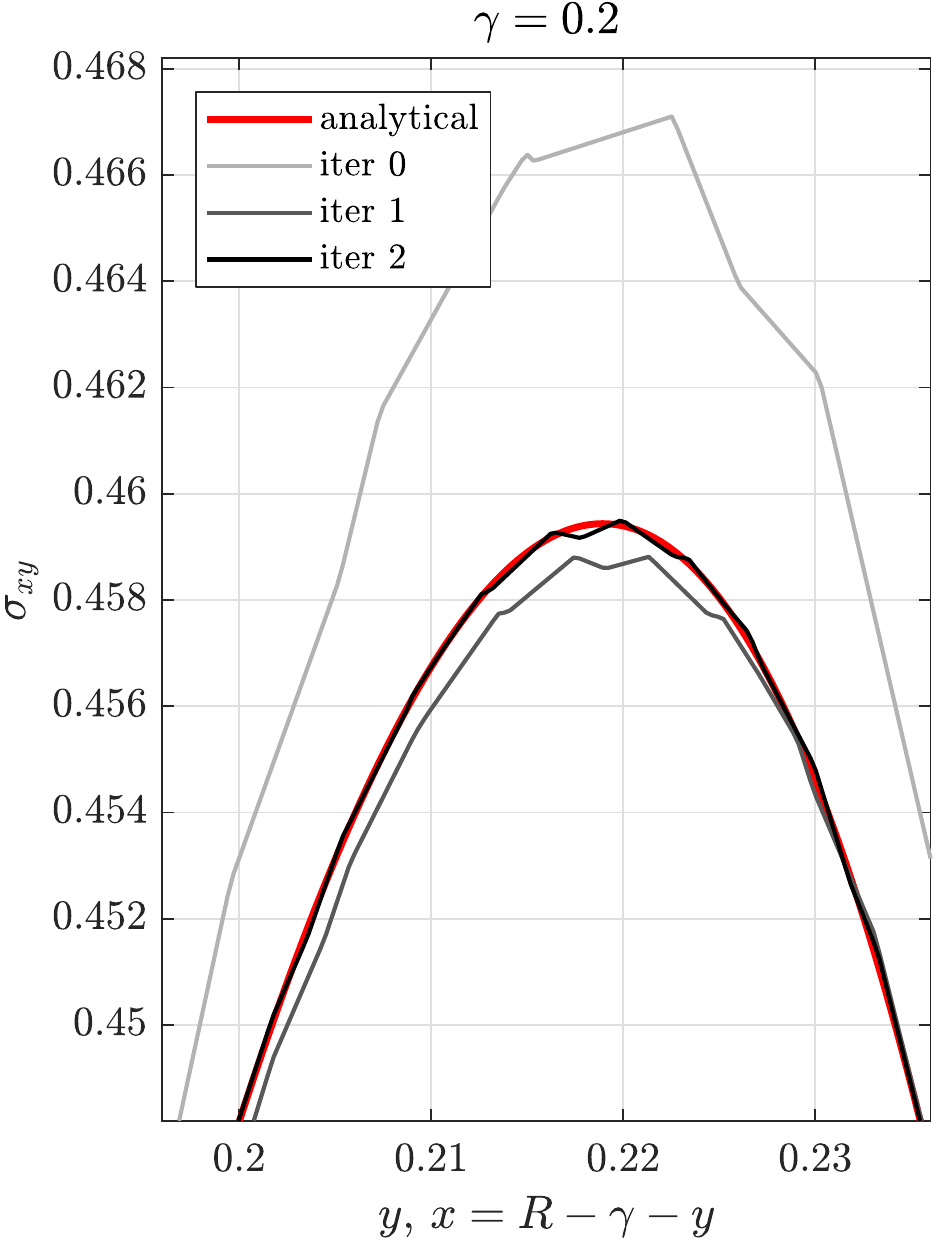}
  \includegraphics[width=0.25\linewidth]{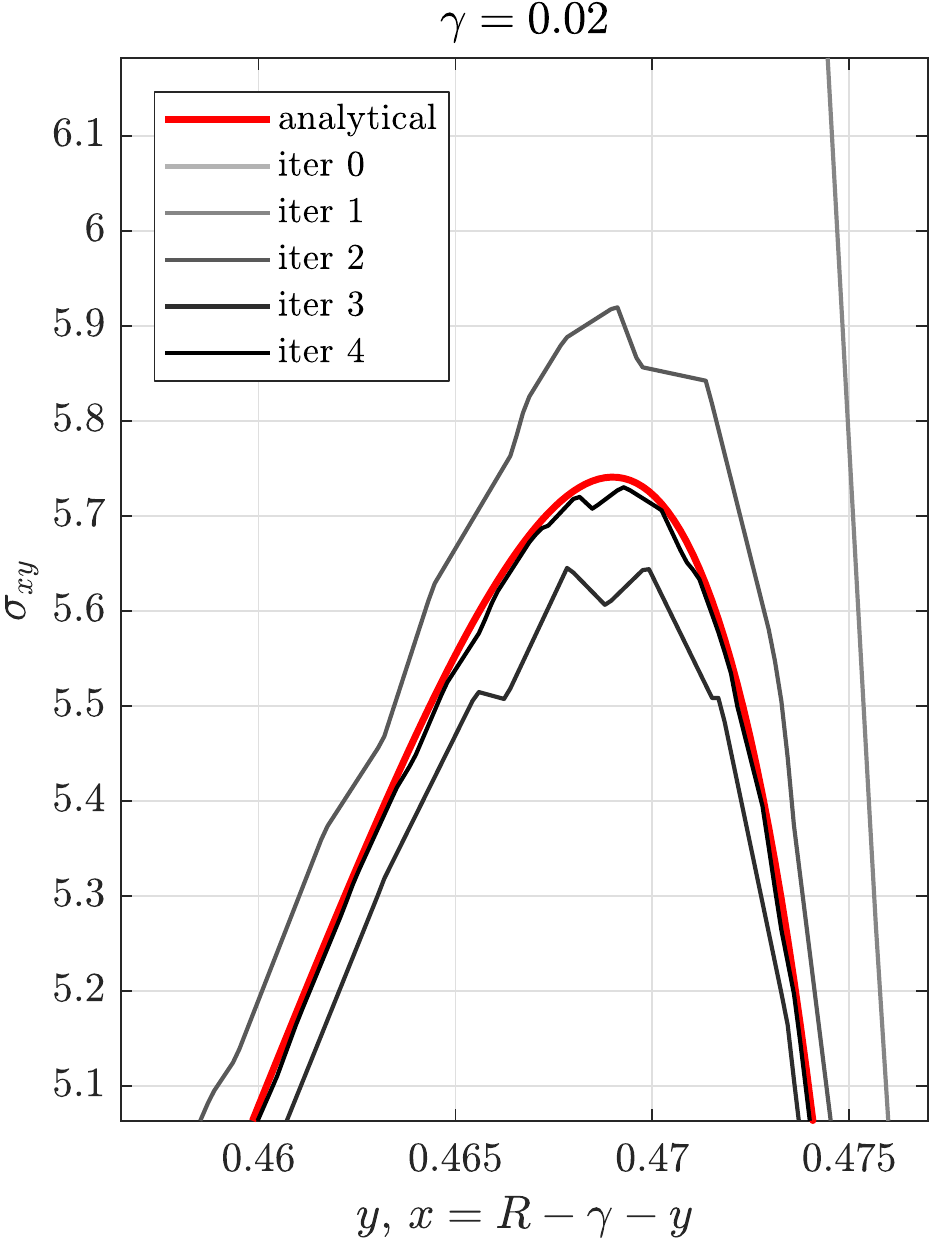}
  \includegraphics[width=0.25\linewidth]{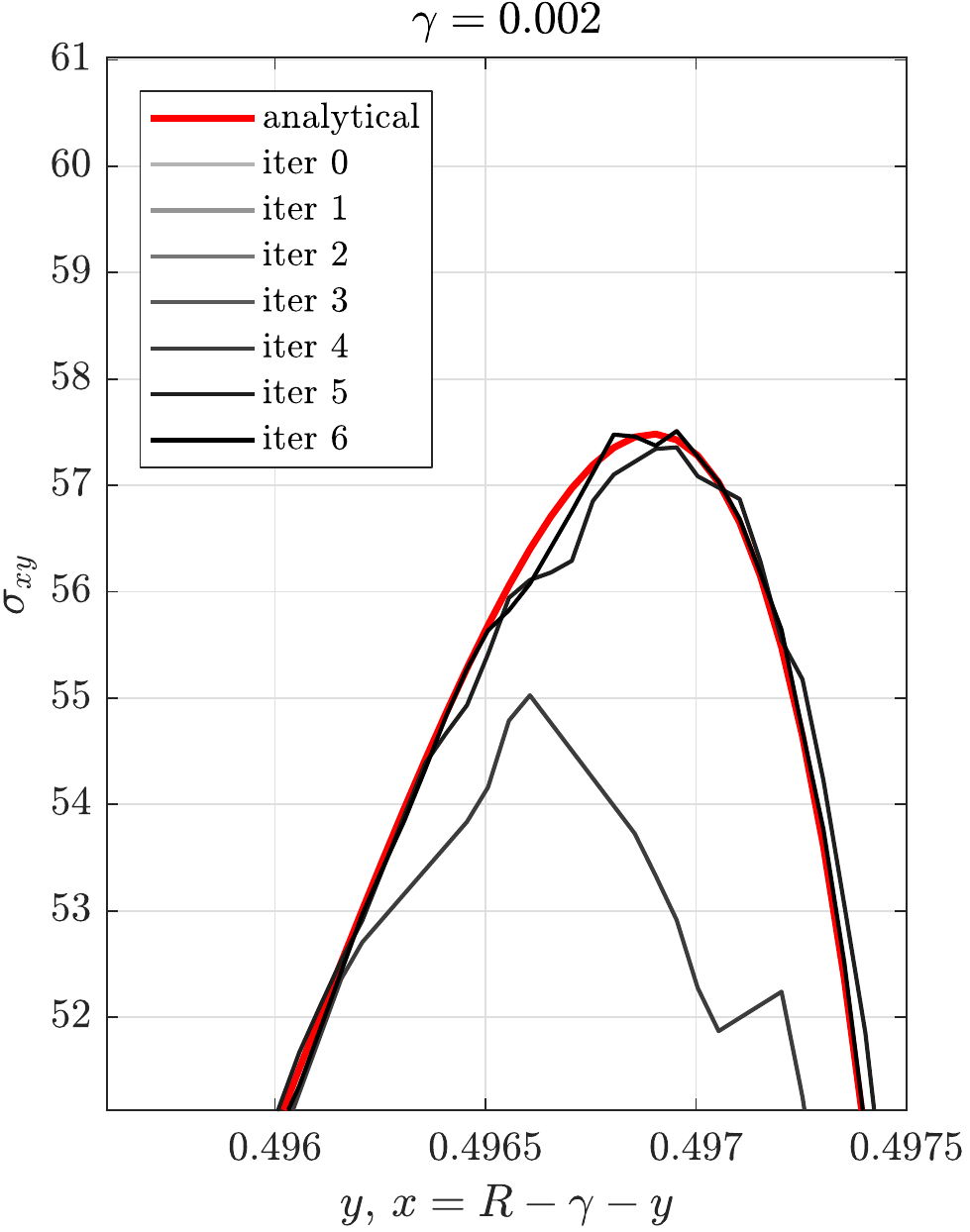}
  \caption{Computed $\sigma_{xy}$ stress profiles for three considered cases. The lower three
    images show close-ups of the computed peaks.}
  \label{fig:profiles}
\end{figure}

The stress profiles computed with the presented adaptive procedure for all three cases are
shown in Fig.~\ref{fig:profiles}. To avoid cluttering, only $\sigma_{xy}$
stress component is shown, but all behave similarly. In the first case, the
initial distribution is almost sufficiently dense and thus adaptive procedure
finishes after 2 iterations. The lower panel still reveals gradual improvement
in both iterations. In the second case, the initial solution underestimates the
peak stress value, the next iteration overestimates it, but then the profiles
start to match with the analytic profile, as seen in the lower panel. This
behaviour is even more prominent in the final case, where the initial
distribution covers high stress area very coarsely and the first few iterations
are only beginning to accurately cover the narrow stress peak. The lower panel once
again shows that in the later iterations, when the solution behaviour has been well captured,
the procedure behaves similarly to the other cases.

\begin{figure}[h!]
  \centering
  \includegraphics[width=0.32\linewidth]{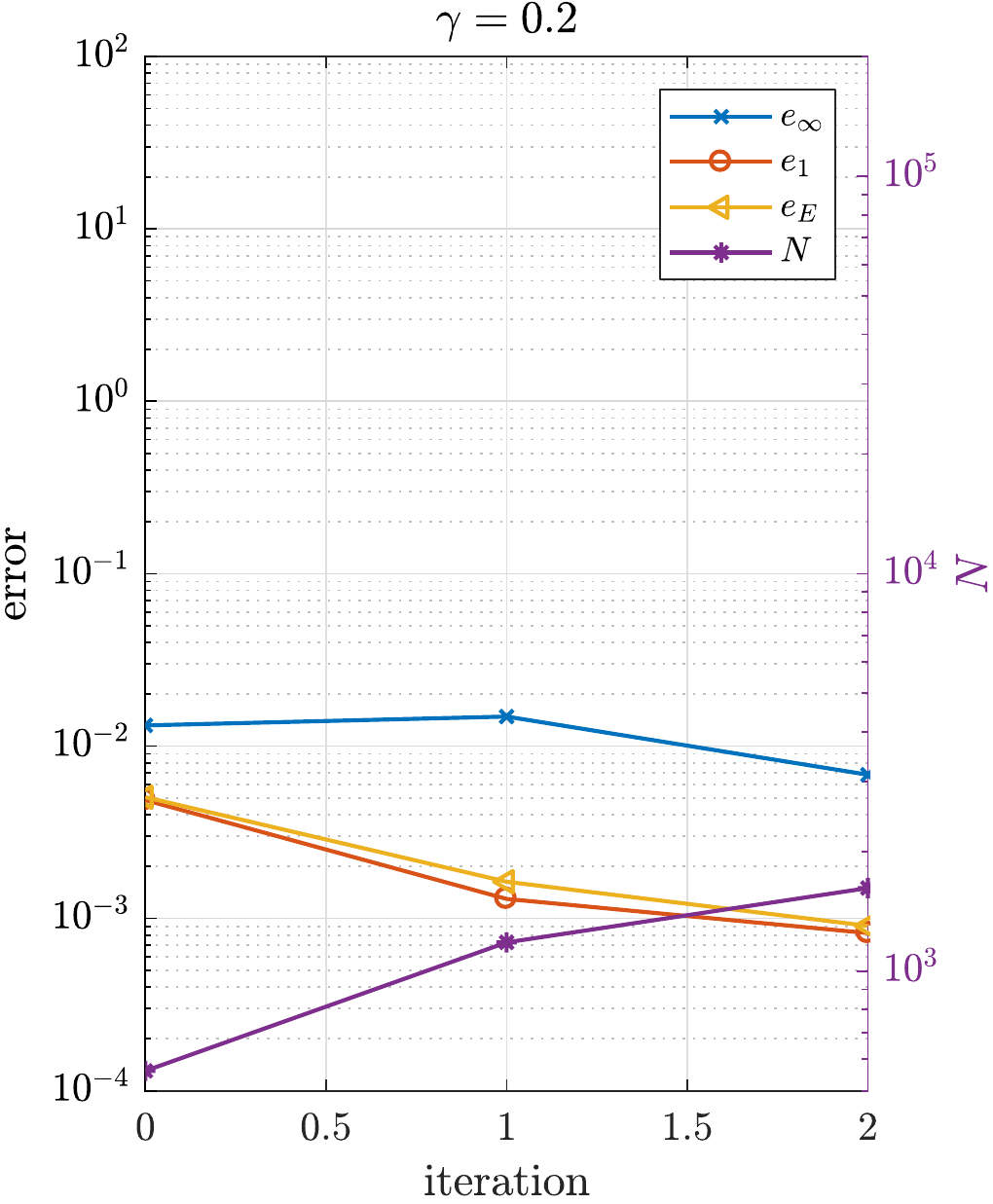}
  \includegraphics[width=0.32\linewidth]{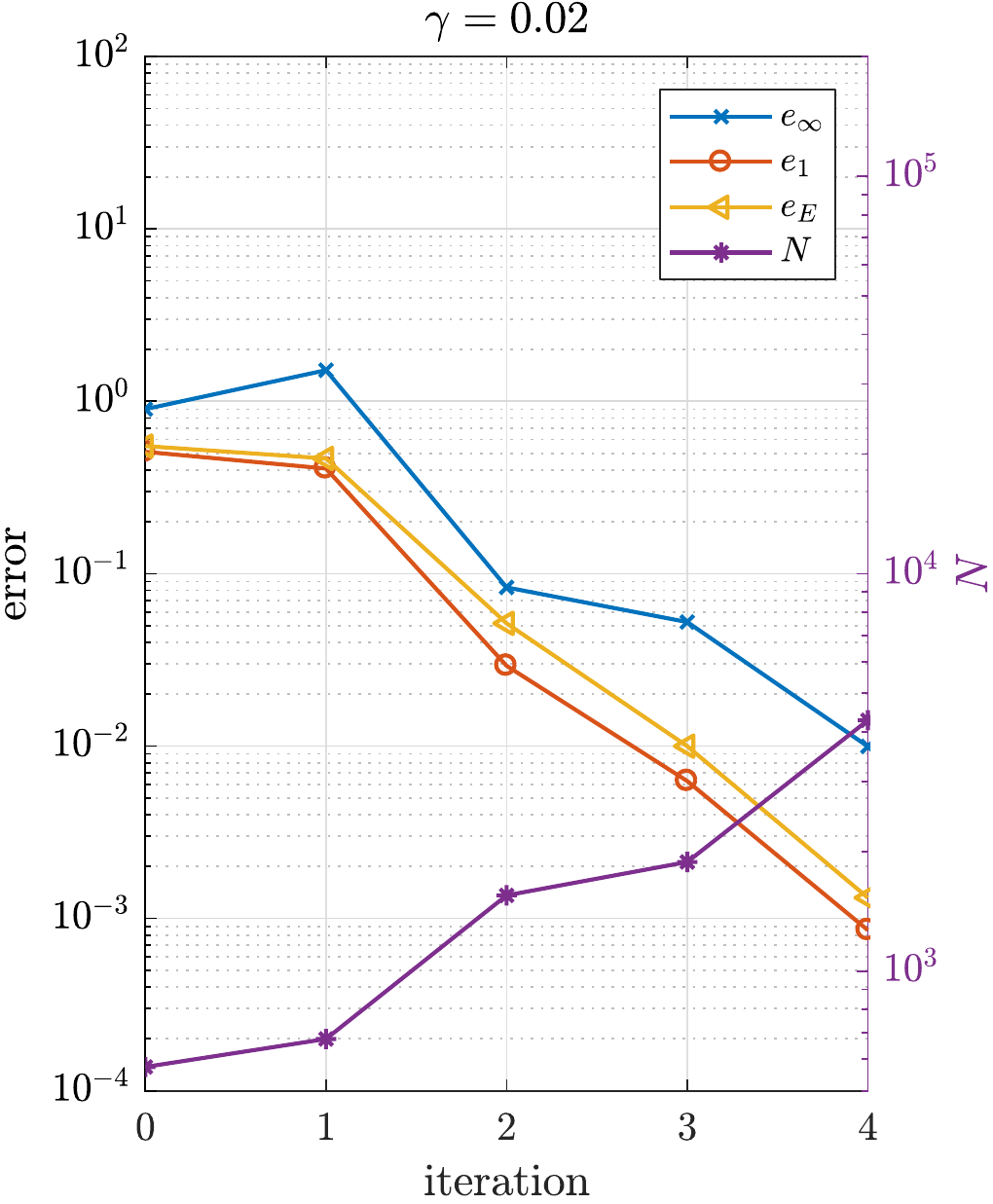}
  \includegraphics[width=0.32\linewidth]{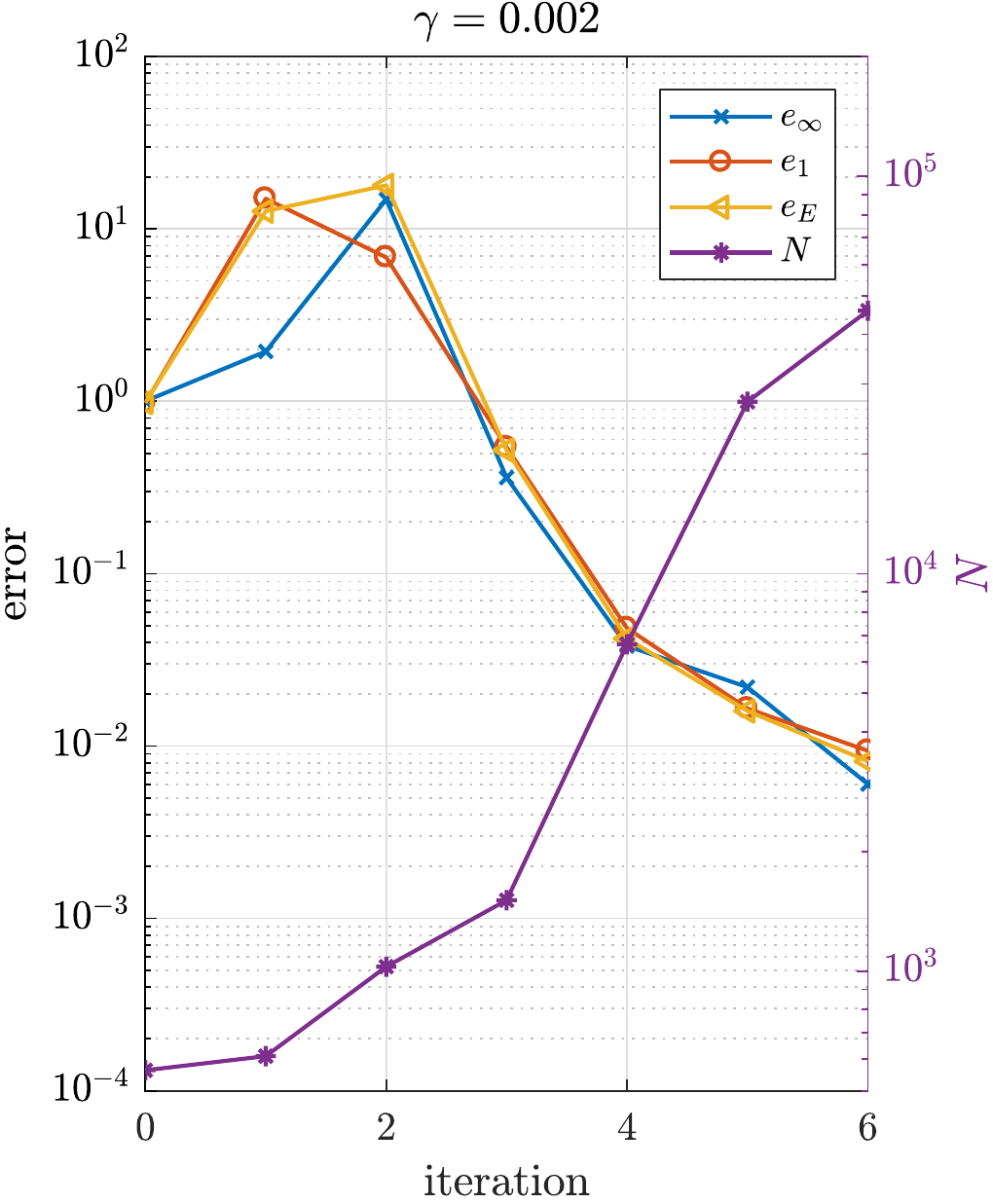}
  \caption{Error of the numerical solution of the compressed disk problem during adaptive iteration.}
  \label{fig:ac}
\end{figure}

The errors over the course of the adaptive iteration are shown in Fig.~\ref{fig:ac}
and confirm the behaviour observed in Fig.~\ref{fig:profiles}.
All three error norms decrease as the algorithm proceeds,
however no guarantee of monotonicity is given and sudden jumps in error can happen.
In the second and more prominently in the third case, an initial rise in error
can be observed, caused by low initial nodal densities. However in few iterations,
the algorithm increases the density enough to accurately detect the error regions.
Note that the error threshold $\eps$ has no direct connection to the final errors
observed, only a general rule that lowering $\eps$ will also decrease the measured errors holds.
The number of nodes is also shown and behaves as expected: mode difficult cases require
more iterations to obtain the final solution, which also uses more nodes for harder cases.

\begin{figure}[h!]
  \centering
  \includegraphics[width=0.8\linewidth]{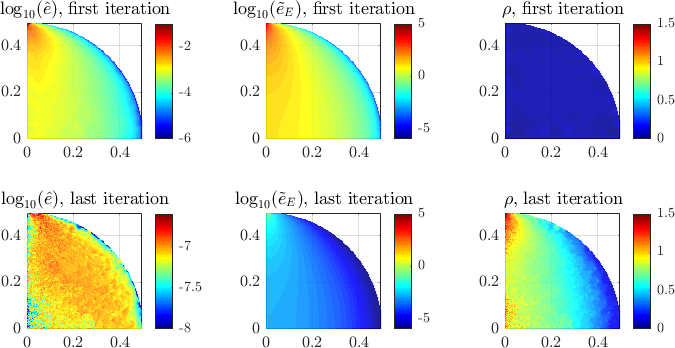}
  \caption{Error indicator value in the last iteration.}
  \label{fig:err-den}
\end{figure}

To better visualise the adaptive results, a comparison between initial
and final values of errors, error indicator and nodal densities are shown in
Fig.~\ref{fig:err-den}.

The measure of relative nodal density $\rho$ at $p_i$ is defined as
\begin{equation}
\rho_i = -\log_{10}(\dr_i/\max_{i}(\dr_i)),
\end{equation}
where $\dr_i$ is the distance to the closest neighbouring node. A value of e.g.~$\rho_i = 2$
means that nodal density around node $p_i$ is $10^2$ times greater than the lowest density
in the domain. For the measure of error at node $p_i$ the kernel of the energy error
$\tilde{e}_E$~\eqref{eq:energy-norm} was taken.
Initially the nodal density is constant (i.e. $\rho \approx 0$) and the error and
error indicator approximately match (up to a factor). In the last iteration,
the error indicator is more uniformly distributed over the whole domain and its value is lower in general,
as expected. When comparing the error and nodal densities the same colour axis limits
were used for the first and the last iteration to get a correct global impression.
Over the course of the adaptive iteration the high error areas were repeatedly refined,
gradually lowering the error to a much more uniform distribution. The final
density distribution however transformed opposite to the error, from a uniform
distribution to the final distribution that resembles the initial error distribution.
This mutually reverse transition to and from uniformity between error
and node density was observed in all three considered cases.

\begin{figure}[h]
  \centering
  \includegraphics[width=0.32\linewidth]{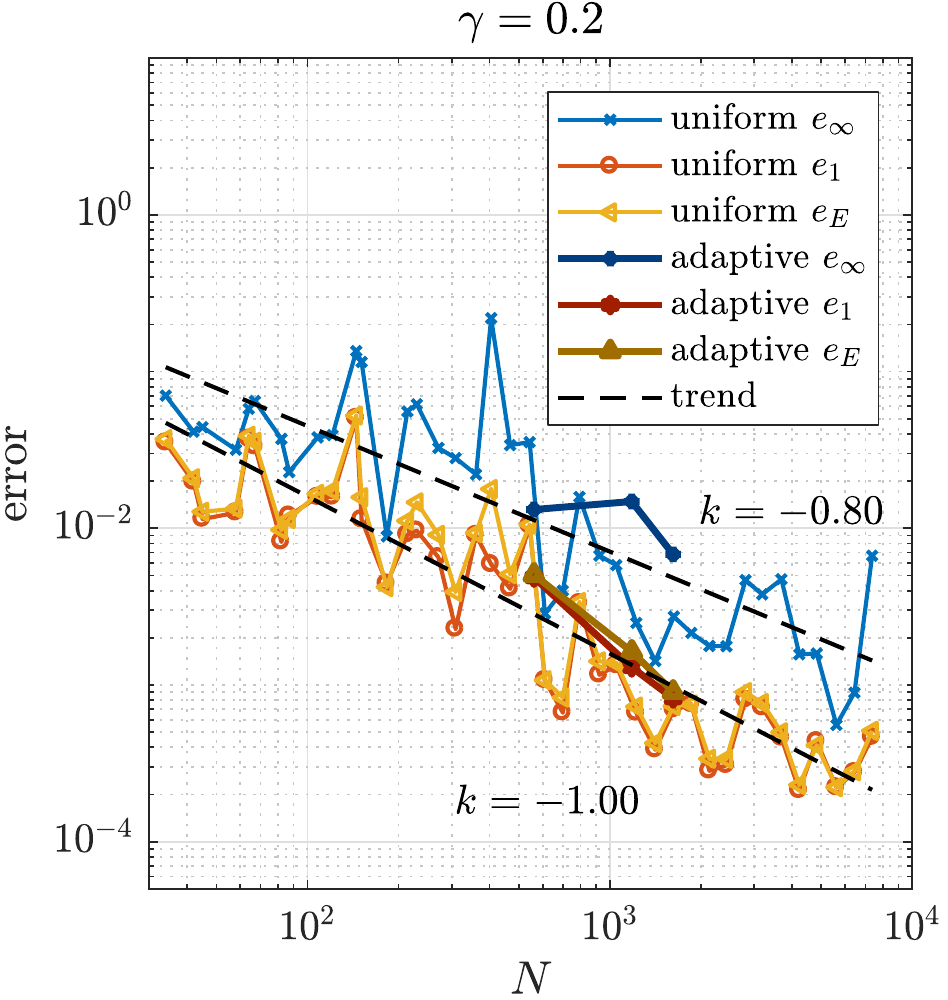}
  \includegraphics[width=0.32\linewidth]{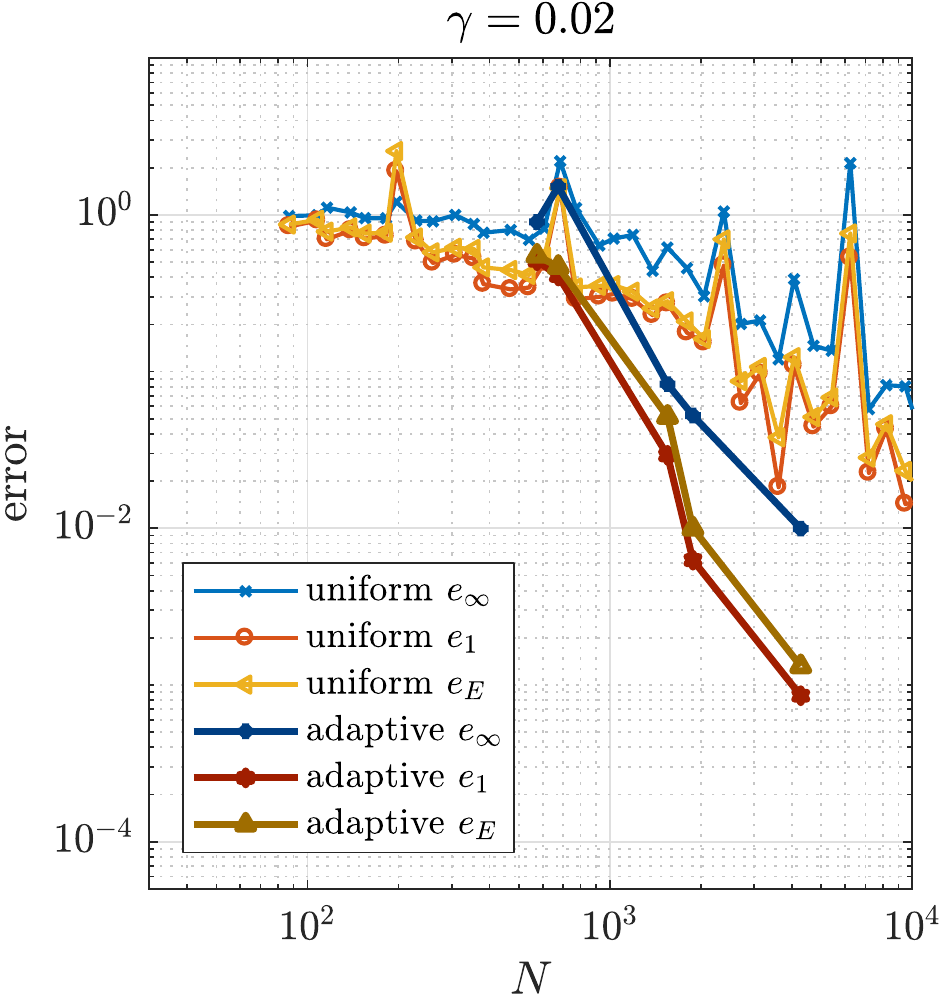}
  \includegraphics[width=0.32\linewidth]{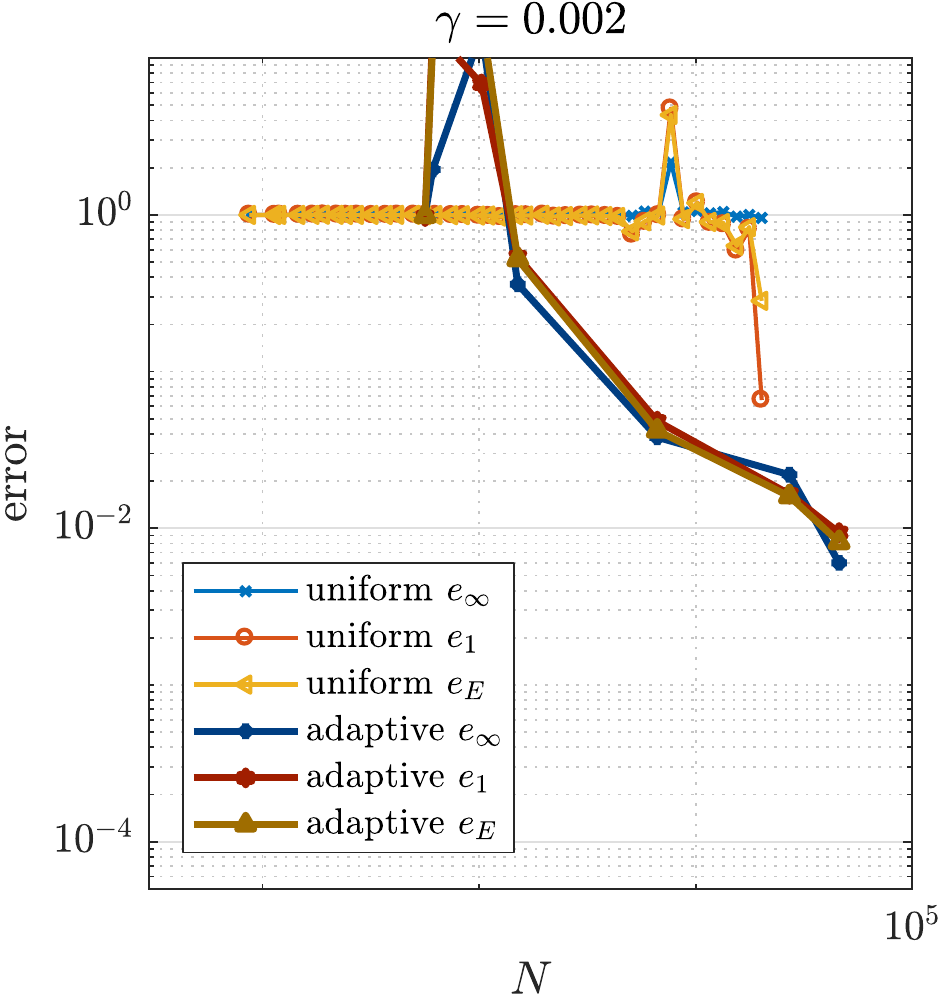}
  \caption{Error of the numerical solution of the compressed disk problem  with respect to number of nodes.}
  \label{fig:uc}
\end{figure}

Finally, to justify testing the adaptivity on these cases at all, the
problem~\eqref{eq:compressed-disk-problem} is solved with uniform node
density, i.e.\ the density function $\dr$ is constant.
Error of the numerical solution was analysed for $\dr$ ranging from
$0.05$ to $0.0033$. The results for all three cases are shown in Fig.~\ref{fig:uc}.

The value of parameter $\gamma$ clearly influences the case difficulty and
the need for adaptivity. In the easiest case uniformly distributed nodes perform well,
and the observed order of convergence in $e_\infty$ matches expectations
provided by the theory for standard finite difference approximations.
The middle case however exhibits worse convergence properties, and the method
completely fails to even compute the most difficult case with satisfiable
accuracy. Adaptivity is clearly needed in the last two cases and
the tested procedure displays satisfactory error behaviour in all three cases.

\subsection{Analysis of adaptivity parameters}
\label{sec:adaptivity-params}
The two important parameters of the refine procedure are the error threshold
$\eps$ and refine aggressiveness $\alpha$. The error threshold parameter $\eps$
is straightforward, and it generally holds that the lower the $\eps$, the lower
the final error will be. Note that if $\eh$ is merely an error indicator and not
an error estimator (i.e.\ it only gives an indication of areas with high errors
and not upper bounds on errors), then the final error may be larger than $\eps$.

The parameter $\alpha$ controls the aggressiveness of the refine procedure, as it
bounds the density increase from above by a factor of at most $\alpha$. Setting
e.g.\ $\alpha = 2$ allows the internodal distance to at most halve and the number
of nodes in a region to increase at most 4 times (in 2D). Intuitively, lower values of
$\alpha$ should result in slower and more controlled error decreases, while higher
values lower the number of needed iterations but the error behaviour is more volatile.

\begin{figure}[h]
  \centering
  \begin{subfigure}[t]{0.32\textwidth}
    \includegraphics[width=0.95\linewidth]{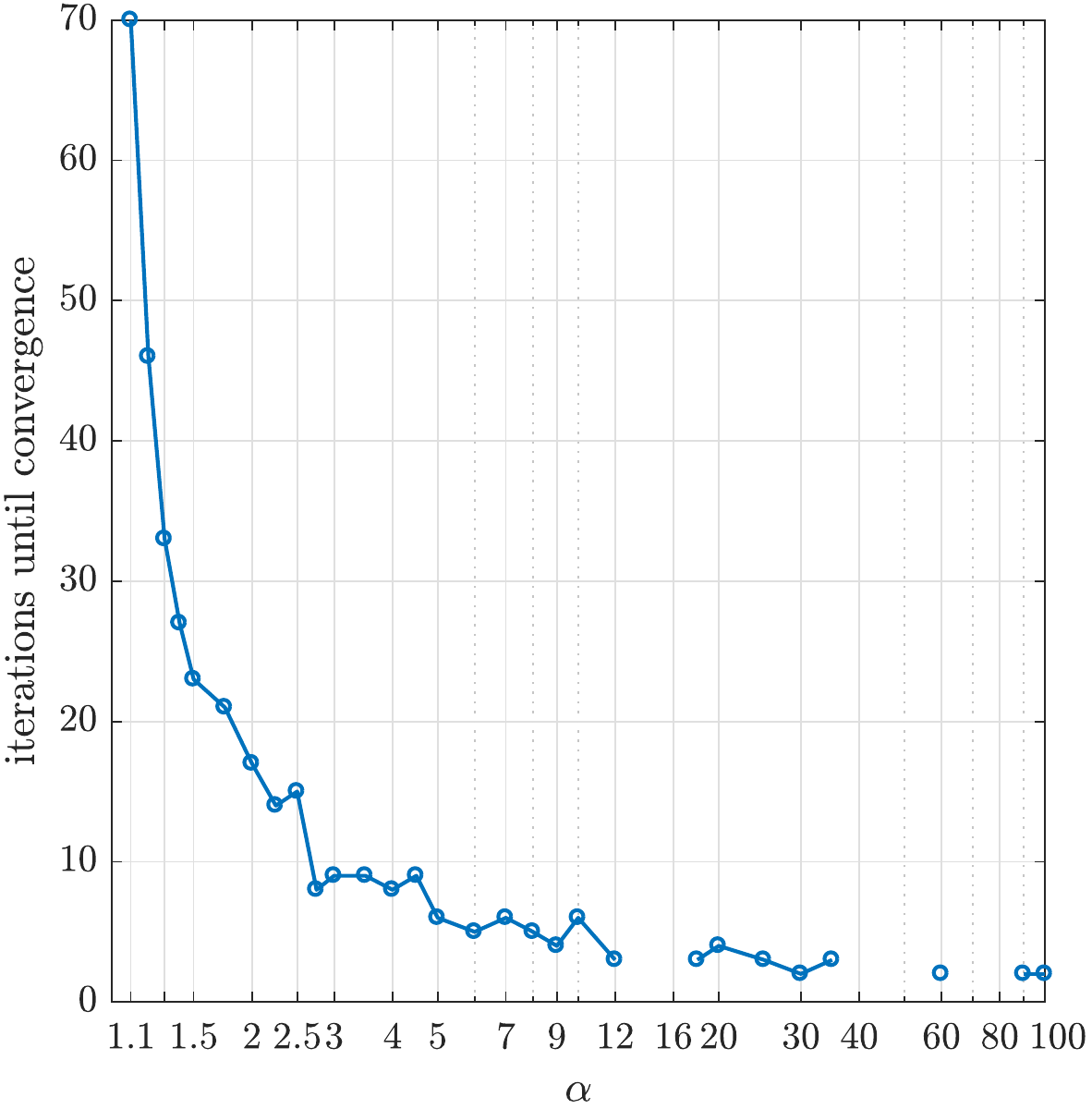}
    \caption{Number of iterations until convergence.}
    \label{fig:alpha-iter}
  \end{subfigure}
  \begin{subfigure}[t]{0.32\textwidth}
    \includegraphics[width=0.95\linewidth]{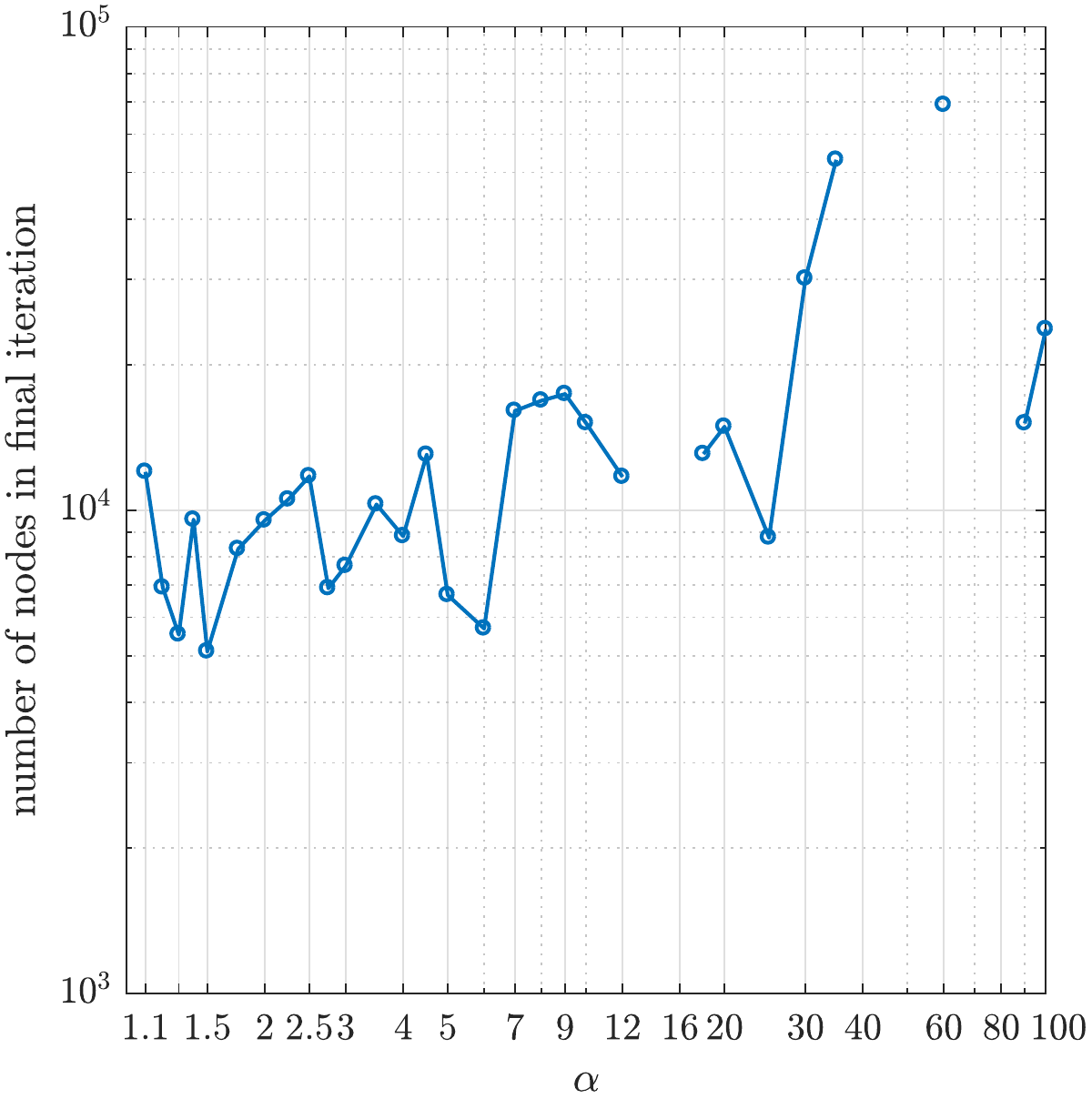}
    \caption{Number of nodes in final iteration.}
    \label{fig:alpha-nodes}
  \end{subfigure}
  \begin{subfigure}[t]{0.32\textwidth}
    \includegraphics[width=0.95\linewidth]{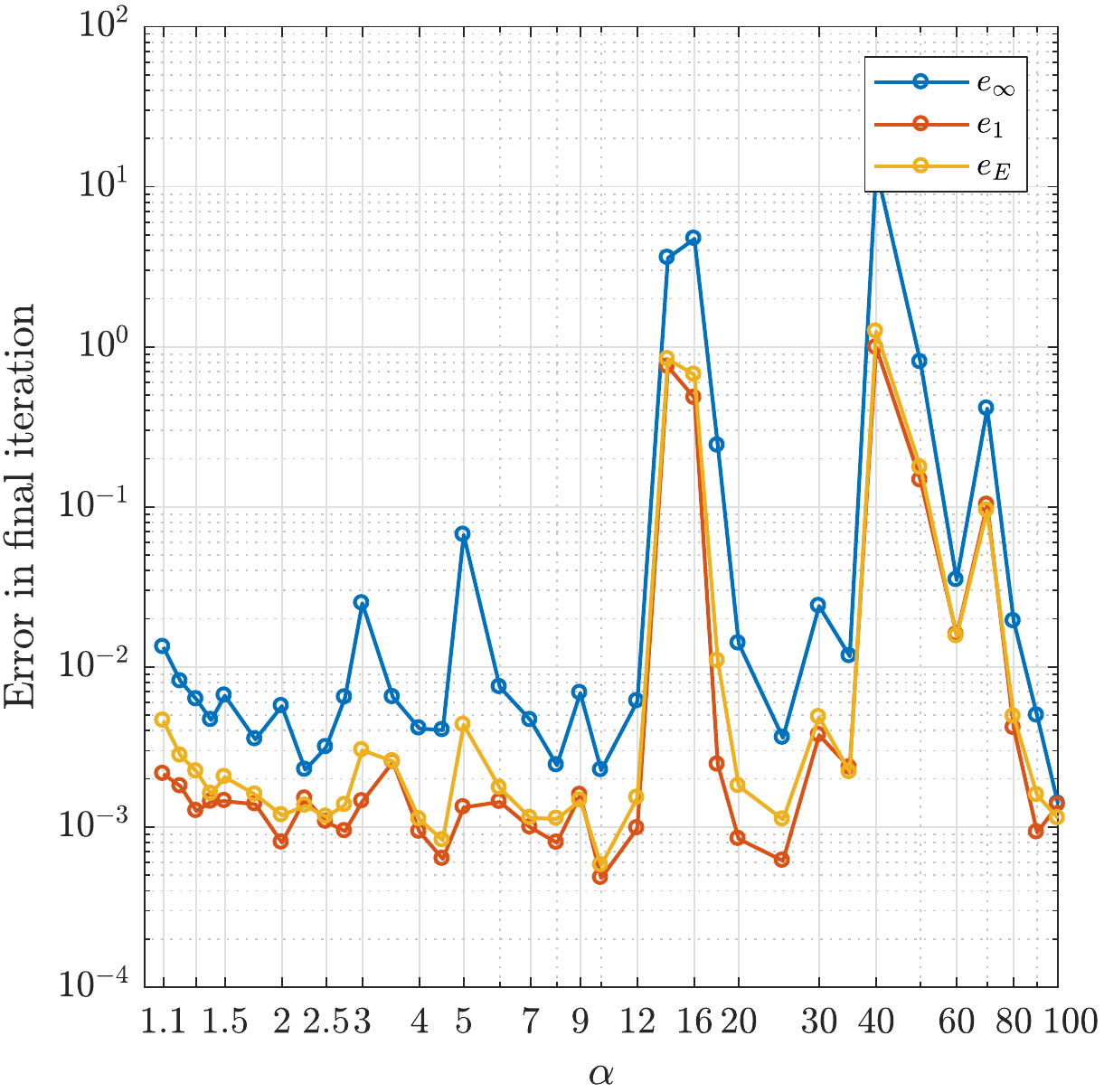}
    \caption{Error in final iteration of each run.}
    \label{fig:alpha-err}
  \end{subfigure}
  \caption{Analyses of effect of aggressiveness parameter $\alpha$.
    Missing data points indicate a divergent solution.}
\end{figure}

To analyse the behaviour with respect to $\alpha$ the case $\gamma = \gamma_2$ was
solved using the same parameters as in section~\ref{sec:half2d}, while
$\alpha$ ranged from 1.1 to 100. Fig.~\ref{fig:alpha-iter} shows the number of
adaptive iterations until convergence, i.e.\ the number of while loop iterations
on line~\ref{alg:while} in algorithm~\ref{alg:adapt}. The missing data points indicate
failure of convergence.

As predicted, the number of iterations reduces steadily with
larger $\alpha$ up to around $\alpha = 10$ when the behaviour becomes too unstable
and failures of convergence occur often. In the stable region with $\alpha \lesssim 10$
both the number of nodes in the final iteration (Fig.~\ref{fig:alpha-nodes}) and
the error in the final iteration (Fig.~\ref{fig:alpha-err}) are mostly
independent of $\alpha$, meaning that irrespective of how many iterations it
took to get to a satisfiable solution, the solution has approximately the same
error and uses approximately the same number of nodes.
Therefore, to reduce the number of iterations necessary while being safely on the
stable side, $\alpha = 5$ was chosen.

The corresponding derefine parameters $\beta$ and $\eta$ behave similarly to
$\alpha$ and $\eps$. Too large values of $\beta$ cause volatile behaviour,
while too low values simply mean that some areas may be covered with too many nodes,
which has little effect on convergence but unnecessarily increases computational costs.
Usually the value of $\beta = 1.5$ was chosen to derefine more conservatively.
The initial distribution $\drI$ was usually chosen for the coarseness bound $\drup$.
This means that derefine cannot reduce the density in the first few iterations,
but can still correct possible increases in wrong areas later on.

The remaining free parameters are related to the density reconstruction
(see sec.~\ref{sec:reconstruct-density}) using modified Sheppard's method.
Varying the number of neighbours $n$ from 3 to 15 had no significant effect
on the solution procedure. Weighting over $n = 7$ nodes was used in most runs.

\subsection{Hertzian contact}
\label{sec:hertz}
The next case considers Hertzian contact between an elastic cylinder and a half
plane, as recently used by Slak and Kosec~\cite{slak2018refined}.
The analytical boundary condition on the top boundary is
given by the known pressure distribution
\begin{equation}
p(x) = \begin{cases}
p_0\sqrt{1-\frac{x^2}{a^2}}, & |x| \leq a \\
0, & \text{otherwise}
\end{cases}, \qquad p_0 = \sqrt{\frac{PE^\ast}{\pi R}}, \quad a = 2\sqrt{\frac{PR}{\pi E^\ast}},
\end{equation}
where $P$ is the pressure force, $R$ is the radius of the cylinder and
the combined Young's modulus $E^\ast$ is given by
$\frac{1}{E^\ast} = \frac{1-\nu_1^2}{E_1} + \frac{1-\nu_2^2}{E_2}$, where $E_1$, $\nu_1$
and $E_2$, $\nu_2$ are the material properties of the cylinder and the half plane, respectively.
The remaining boundary conditions are no-displacement conditions at infinity.

\begin{figure}[h]
  \centering
  \includegraphics[width=0.5\linewidth]{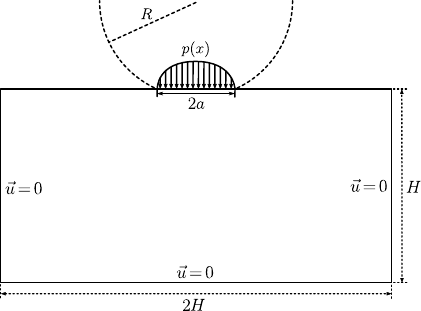}
  \caption{Numerical domain and boundary conditions for the solution of Hertzian contact problem.}
  \label{fig:hertz-bc}
\end{figure}

The problem is solved numerically by truncating the half plane to a rectangle
$[-H, 0] \times [-H, H]$ and applying the no displacement boundary conditions
at all boundaries but the top, as illustrated in Fig.~\ref{fig:hertz-bc}.
For the choice $R = \unit[1]{m}$, $P = \unitfrac[543]{N}{m}$,
$E = \unit[72.1]{GPa}$, $\nu = 0.33$, the contact width $a$ is approximately
$\unit[0.13]{mm}$. A choice of $H = \unit[0.1]{m}$, which is approximately
$770$ times larger than the phenomenon of interest, is sufficiently large
that the truncation error does not present a significant contribution~\cite{slak2018refined}.

The large ratio between domain size and contact width presents a challenging case
for adaptive refinement algorithms. The domain was initially filled with a
constant discretisation step $h$, using the filling algorithm described in
section~\ref{sec:fill}. The initial step $h$ needs to
be small enough to capture the contact itself, while this same $h$ might be
unnecessarily small for the rest of the domain. This offers
a chance to demonstrate adaptive derefinement as well.

Additionally, instead of an analytical indicator given by~\eqref{eq:anal-err-est},
a simple ad-hoc error indicator was tested, to develop a functional
indicator for cases where no analytical solution is available.

The construction of the indicator starts on the premise that the
discretisation needs to be refined in areas with high gradients, i.e.~the areas
where the variability of the numerical solution is locally large. Similar indicators have been already reported for adaptive solution of Poisson and Burgers's equation~\cite{Davydov,kos2}. The error indicator is thus introduced as
\begin{equation}
\hat{e}_i = \underset{p_j \in N_i}{\operatorname{std}} \sigma_{xx}(p_j) +
\underset{p_j \in N_i}{\operatorname{std}} \sigma_{yy}(p_j) +
\underset{p_j \in N_i}{\operatorname{std}} \sigma_{xy}(p_j),
\label{eq:err-indicator-std}
\end{equation}
.

\begin{figure}[h]
  \centering
  \includegraphics[width=0.45\linewidth]{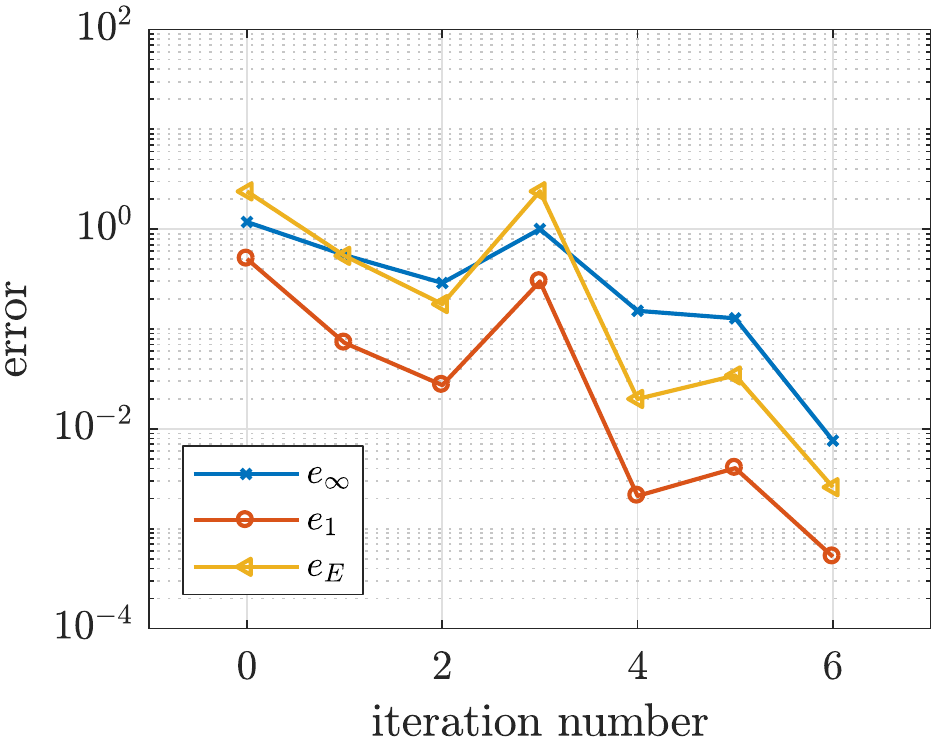}
  \includegraphics[width=0.45\linewidth]{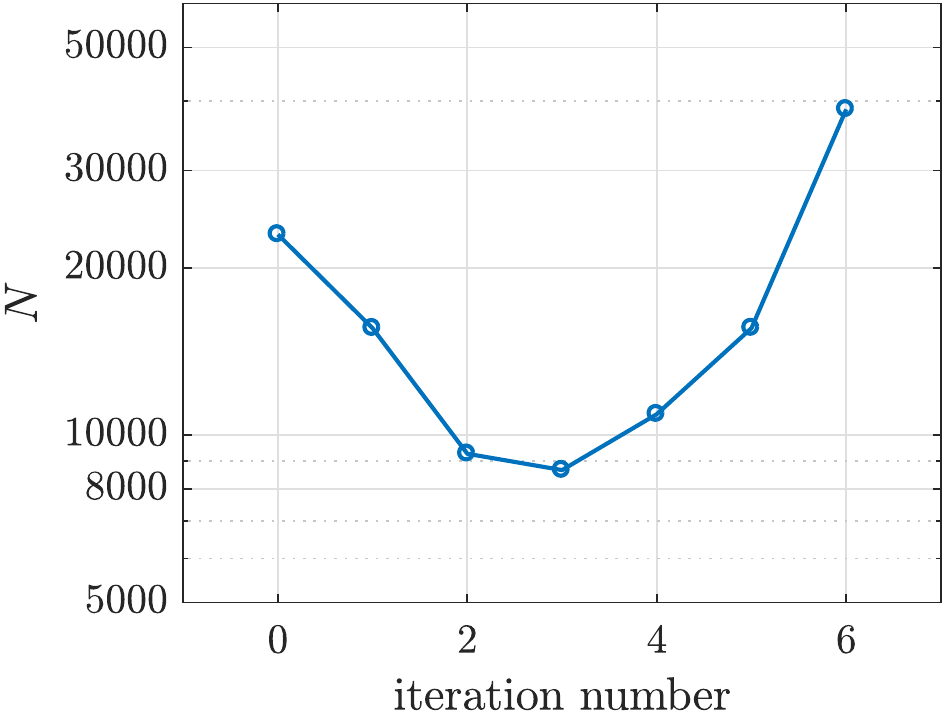}
  \caption{Errors and node counts during adaptive iteration for the solution of the Hertzian problem.}
  \label{fig:hertz-err}
\end{figure}

The adaptive procedure was run with
$\alpha = 5$, $\eps = 10^5$, $\beta = 1.5$, $\eta = 10^4$ and left to run for 6 iterations.
The errors measured against a closed form solution~\cite[eq.~(54--56)]{slak2018refined}
and the number of nodes over the course of the iterative procedure are
shown in Fig.~\ref{fig:hertz-err}. Additionally, data about node refinement
for each iteration is shown in Table~\ref{tab:node-cnt}.

\begin{table}[h]
  \centering
  \begin{tabular}{r|r|r|r|r|l}
    iteration & total nodes & refined & no change & derefined & derefined, but hit bound $\drup$ \\ \hline
    0 & 23\,012 & 384   & 3\,427 & 19\,201 & 0    \\
    1 & 15\,596 & 483   & 777  & 14\,336 & 0   \\
    2 & 9\,264  & 1\,127  & 1\,469 & 6\,668  & 0    \\
    3 & 8\,658  & 5\,198  & 3\,460 & 0     & 0    \\
    4 & 10\,918 & 3\,097  & 2\,017 & 5\,804  & 0    \\
    5 & 15\,605 & 9\,494  & 2\,189 & 3\,753  & 169  \\
    6 & 38\,674 & / & / & / & /  \\
  \end{tabular}
  \caption{Number of refined and derefined nodes for each iteration.}
  \label{tab:node-cnt}
\end{table}

\begin{figure}[h]
  \centering
  \includegraphics[width=0.4\linewidth]{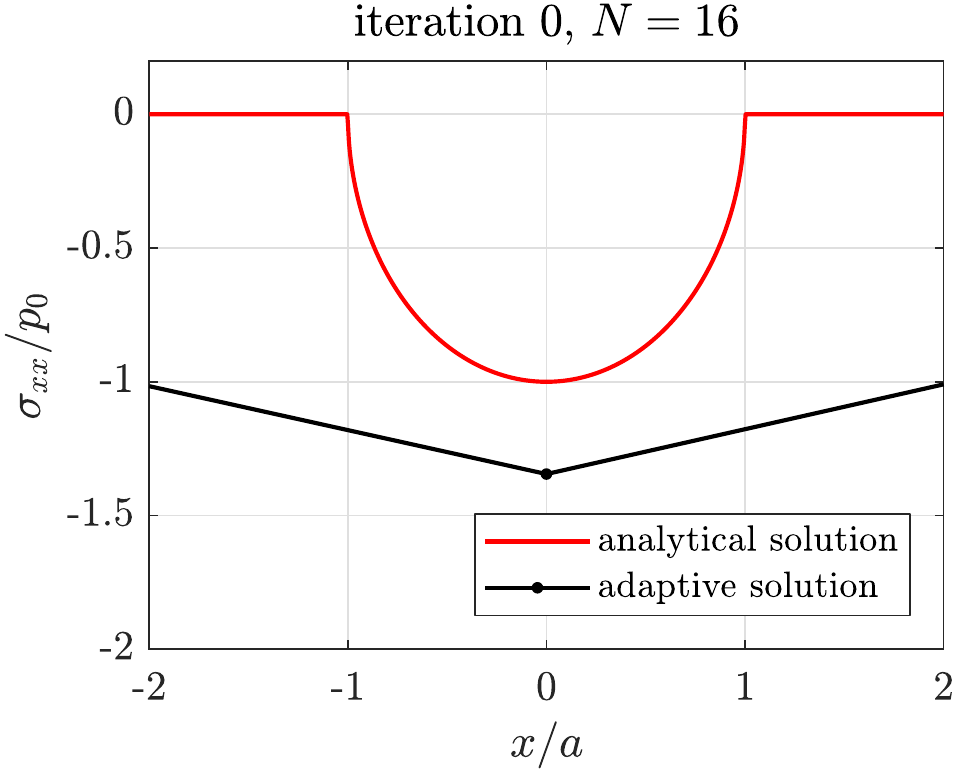}
  \includegraphics[width=0.4\linewidth]{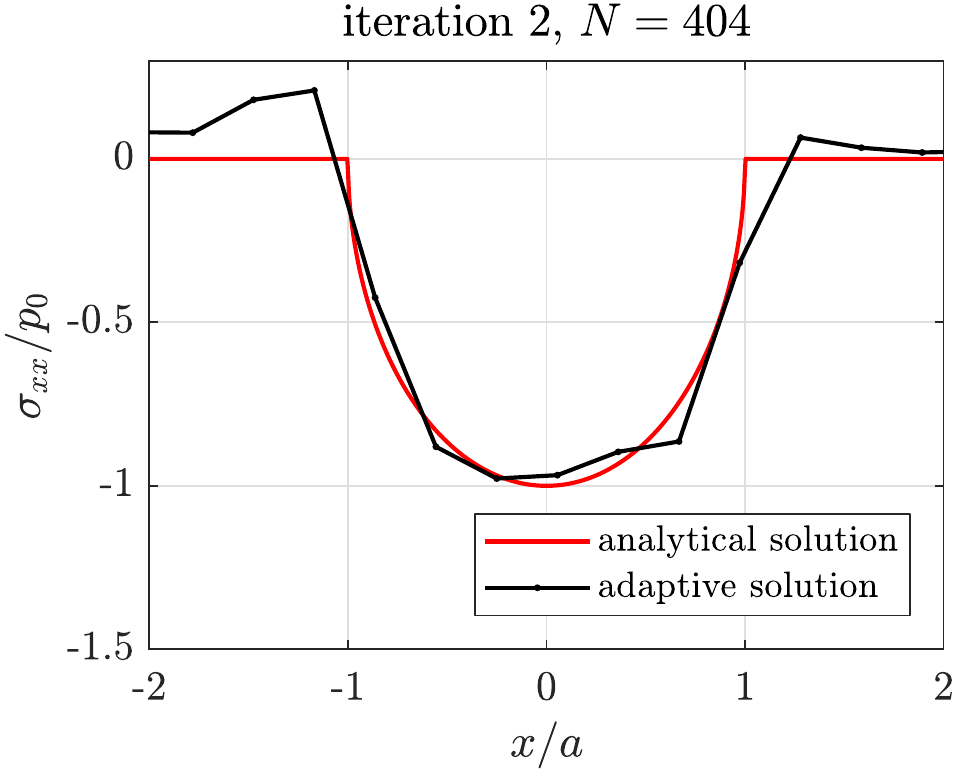}

  \includegraphics[width=0.4\linewidth]{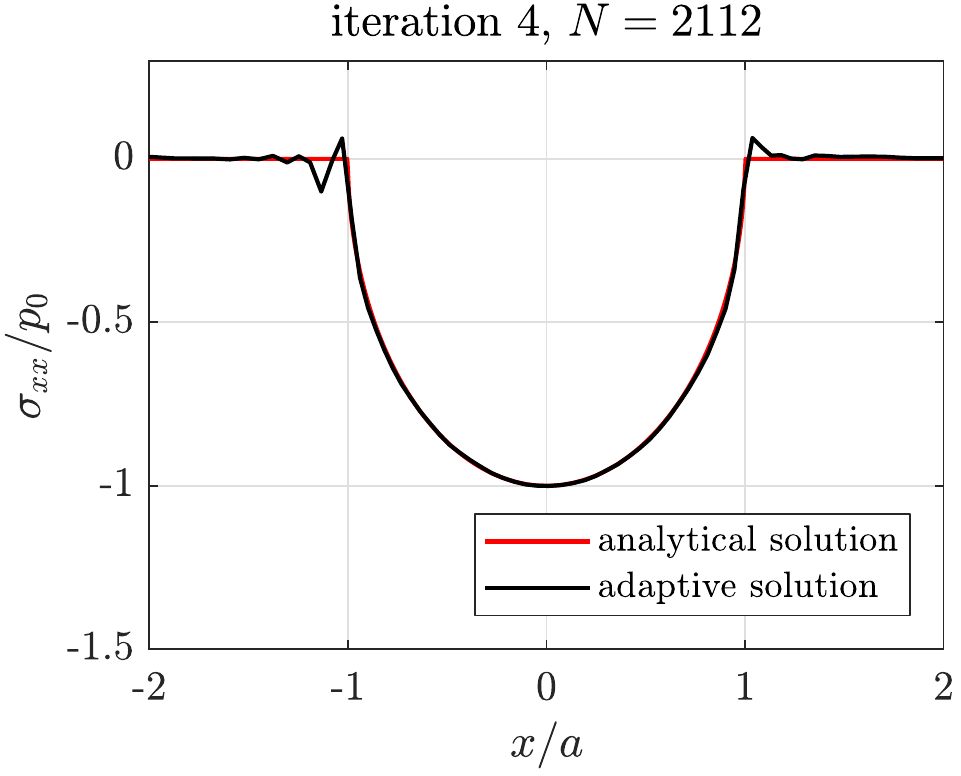}
  \includegraphics[width=0.4\linewidth]{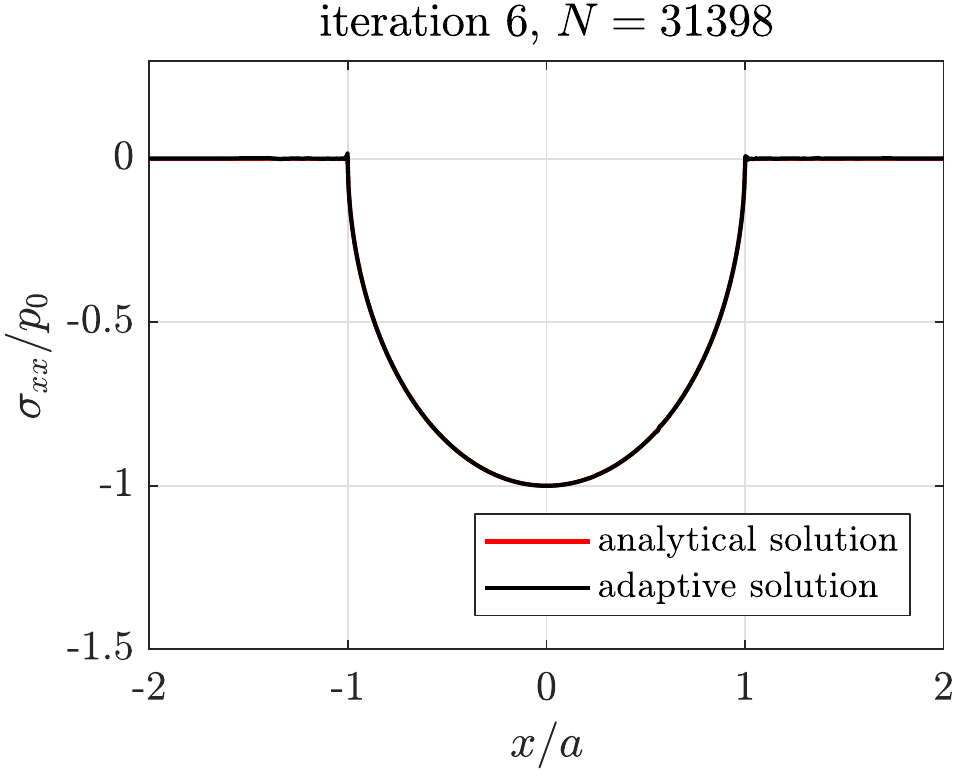}
  \caption{Top surface stress profiles during the iteration. Only even iterations are shown for brevity.}
  \label{fig:hertz-top}
\end{figure}

In the beginning, the derefinement is substantial and the node count decreases,
however, refinement is already present for nodes under contact.
Later in the procedure, refinement becomes more pronounced and the total
number of nodes increases. The computed top surface tractions over
the course of the adaptive iteration are shown in Fig.~\ref{fig:hertz-top}.

\begin{figure}[h!]
  \centering
  \includegraphics[width=0.45\linewidth]{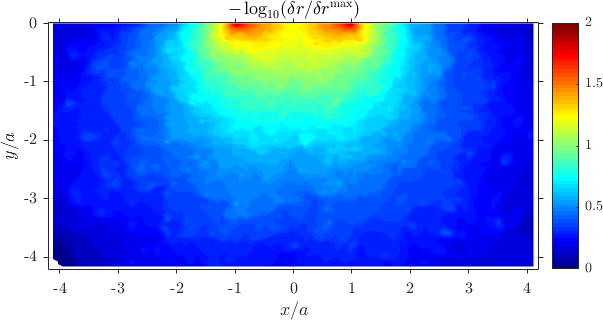}
  \includegraphics[width=0.45\linewidth]{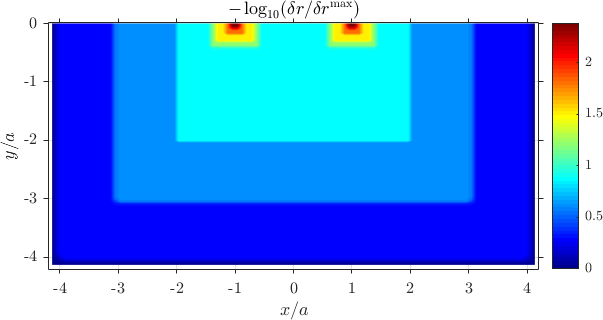}
  \caption{Comparing adaptive (left) and manual (right) nodal distributions densities for the Hertzian problem.}
  \label{fig:hertz-densty}
\end{figure}

It is interesting to compare these results to the manually refined solution obtained
by Slak and Kosec~\cite{slak2018refined}. The error in the final iteration of the adaptive
procedure is in the same range as the errors obtained by manual refinement.
Additionally, an interesting comparison of node densities under contact can be
made and is shown in Fig.~\ref{fig:hertz-densty}.  The manual and adaptively
obtained densities generally refine the same areas, however the
adaptively obtained density is much smoother.

From the error plot in Fig.~\ref{fig:hertz-err} and the density plot in
Fig.~\ref{fig:hertz-densty} it can be concluded that the error indicator
$\eh$, given by~\eqref{eq:err-indicator-std}, works adequately in practice.



\subsection{Fretting fatigue contact}
\label{sec:fwo-case}
\begin{figure}[h]
  \centering
  \begin{subfigure}[t]{0.47\textwidth}
    \centering
    \includegraphics[width=0.95\linewidth]{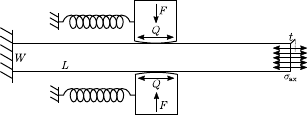}
    \caption{Schema of the experiment.}
    \label{fig:fwo-case-def}
  \end{subfigure}
  \begin{subfigure}[t]{0.47\textwidth}
    \centering
    \includegraphics[width=0.95\linewidth]{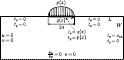}
    \caption{Numerical domain and boundary conditions.}
    \label{fig:fwo-case-bc}
  \end{subfigure}
  \caption{Description of the considered fretting fatigue case. Ratios in the drawings are not to scale.}
  \label{fig:fwo-drawings}
\end{figure}

The usability of the described adaptive technique is further demonstrated on a
case, originating from the study of fretting fatigue, described
by Pereira et al.~\cite{pereira2016on}, for which the authors believe no closed form solution
is known.  A small thin rectangular specimen of width $W$, length $L$ and
thickness $t$ is stretched in one axis with axial traction $\sigma_{ax}$ and
compressed in another with force $F$ by two oscillating cylindrical pads of
radius $R$, inducing tangential force $Q$,  as shown in
Fig.~\ref{fig:fwo-case-def}. The tractions induced by the pads are predicted
using an extension of Hertzian contact theory, which splits the contact area into
the stick and slip zones, depending on the coefficient of friction $\mu$ and
combined Young's modulus $E^*$, given by $\frac{1}{E^*} = \frac{1 -
  \nu_1^2}{E_1} + \frac{1 - \nu_2^2}{E_2}$, where $E_i$ and $\nu_i$ represent the
Young's moduli and Poisson's ratios of the specimen and the pad, respectively.
The theory predicts the contact half-width
\begin{align}
a = 2 \sqrt{\frac{FR}{t\pi E^*}},
\end{align}
normal traction
\begin{equation}
p(x) = \begin{cases}
p_0 \sqrt{1-\frac{x^2}{a^2}}, & |x| < a \\
0, & \text{else}
\end{cases}, \qquad p_0 = \sqrt{\frac{FE^*}{t \pi R}},
\end{equation}
and tangential traction
\begin{equation}
q(x) = \begin{cases}
-\mu p_0 \left(\sqrt{1 - \frac{x^2}{a^2}} - \frac{c}{a}\sqrt{1 - \frac{(x-e)^2}{c^2}}\right), &
|x-e| < c \\
-\mu p_0 \sqrt{1 - \frac{x^2}{a^2}}, & c \leq | x - e | \text{ and } |x| \leq a \\
0, & \text{else}
\end{cases}
\end{equation}
where $c = a\sqrt{1 - \frac{Q}{\mu f}}$ is the half-width of the slip zone, and
$e = \operatorname{sgn}(Q)\frac{a \sigma_{ax}}{4 \mu p_0}$ is the
eccentricity due to axial loading. Note that the inequalities $Q \leq \mu F$ and
$\sigma_{ax} \leq 4\left(1 - \sqrt{1 - \frac{Q}{\mu F}}\right)$ must hold
for these expressions to be valid, both of which will be satisfied in our case.

Plane strain is used to reduce the problem to two dimensions and symmetry along the
horizontal axis is used to reduce the problem size. The region
$\Omega = [-L/2, L/2] \times [-W/2, 0]$ is taken as the problem domain and the
boundary conditions are illustrated in Fig.~\ref{fig:fwo-case-bc}.

\begin{figure}[h]
  \centering
  \includegraphics[width=0.8\textwidth]{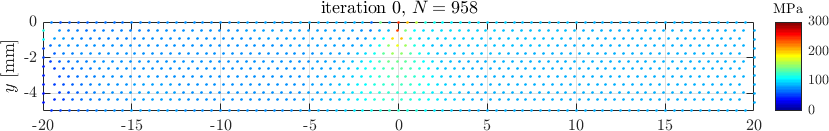} \\[1ex]
  \includegraphics[width=0.8\textwidth]{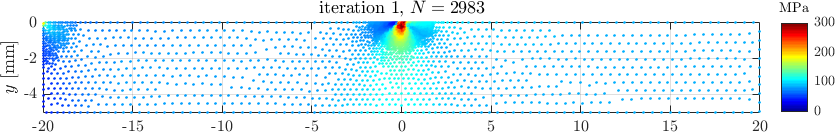} \\[1ex]
  \includegraphics[width=0.8\textwidth]{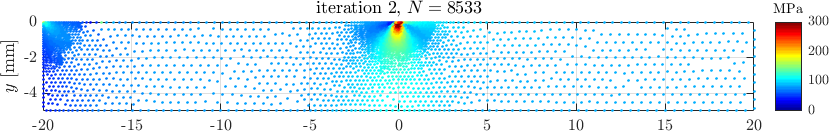} \\[1ex]
  \includegraphics[width=0.8\textwidth]{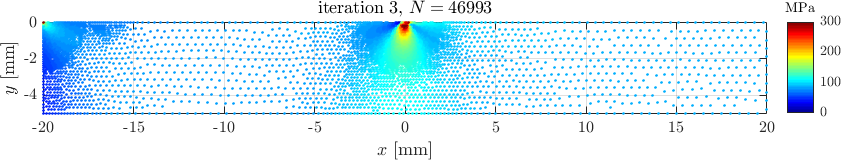}
  \caption{Nodes of four subsequent adaptive iterations during the solution of the fretting fatigue problem, coloured by von Mises stress.}
  \label{fig:fwo-adapt-mesh}
\end{figure}

Values of $E_1 = E_2 = \unit[72.1]{GPa}$, $\nu_1 = \nu_2 = 0.33$ were taken
for the material parameters, coinciding with aluminium 2420-T3. Dimensions
of the specimen were $L = \unit[40]{mm}$, $W = \unit[10]{mm}$ and $t = \unit[4]{mm}$.
Values $F = \unit[543]{N}$, $Q = \unit[155]{N}$, $\sigma_{ax} = \unit[100]{MPa}$,
$R = \unit[10]{mm}$ and $\mu = 0.3$ were chosen for the model parameters.
This means that the half-contact width $a$ equals $\unit[0.2067]{mm}$, which is
approximately 200 times smaller than the domain width $W$, once again making it
difficult to solve without adaptivity.
The same ad-hoc error indicator as in section~\ref{sec:hertz} was used to adaptively
refine the solution. The initial distribution was generated using constant spacing
$\drI = \unit[0.5]{mm}$, putting exactly one node on the whole contact surface.
4 iterations of the adaptive solution procedure were run
using $\eps = \unit[0.1]{MPa}$ as the error threshold and
$\alpha = 10$, resulting in a domain populated with 49\,993 nodes.
The sets of nodes $\P$ and the corresponding numerical solutions
obtained during the adaptive iteration are shown in Fig.~\ref{fig:fwo-adapt-mesh}.
The results are compared
to a solution obtained with the Finite Element Method (FEM)
on a much denser mesh with more than 100\,000 DOFs,
using commercial ABAQUS\textsuperscript{\textregistered}
software for finite element analysis, where manual refinement has been used~\cite{pereira2016on},
as well as to another FEM solution, obtained using the freely available
FreeFem++~\cite{hecht2012ff} and its built-in adaptive refinement techniques.

\begin{figure}[h]
  \centering
  \includegraphics[width=0.3\textwidth]{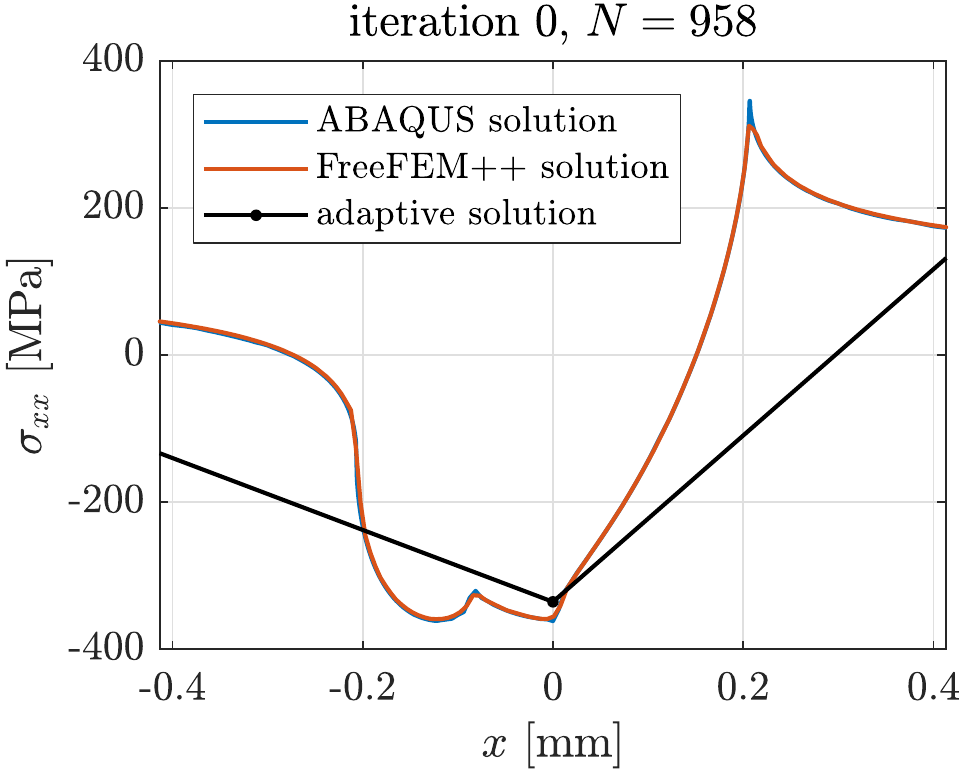}
  \includegraphics[width=0.3\textwidth]{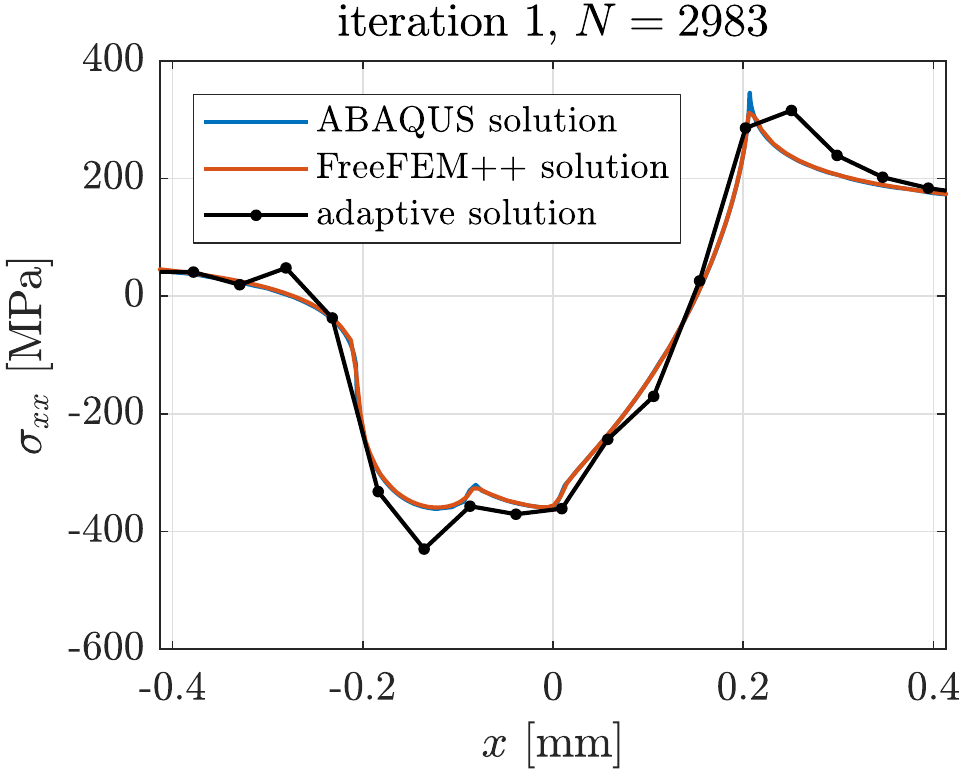} \\[1ex]
  \includegraphics[width=0.3\textwidth]{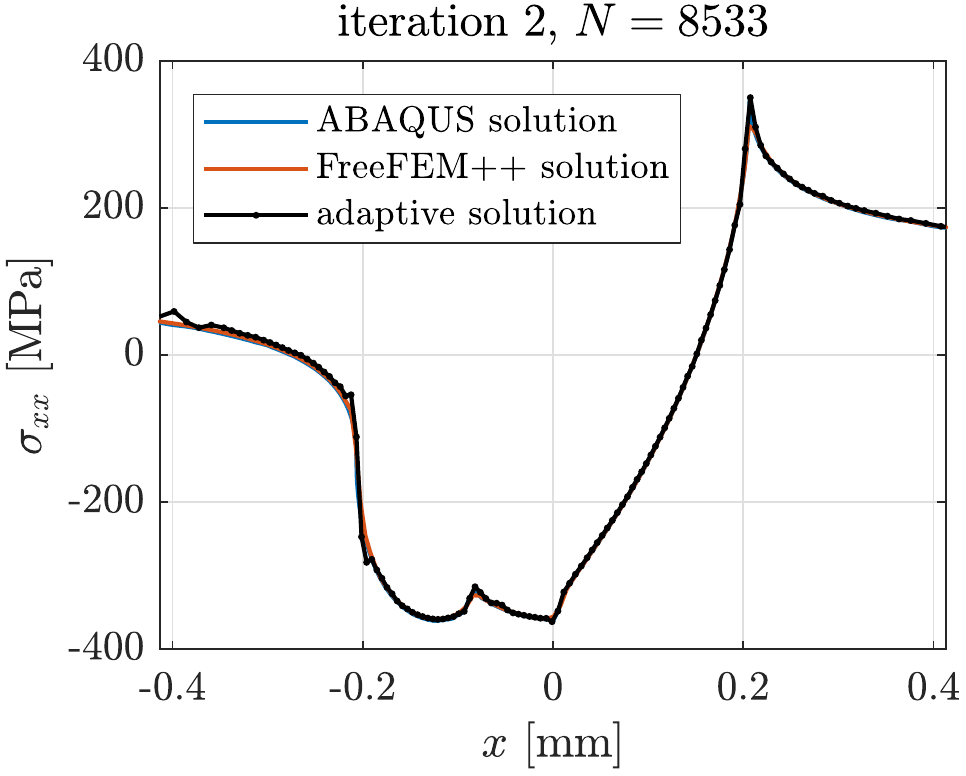}
  \includegraphics[width=0.3\textwidth]{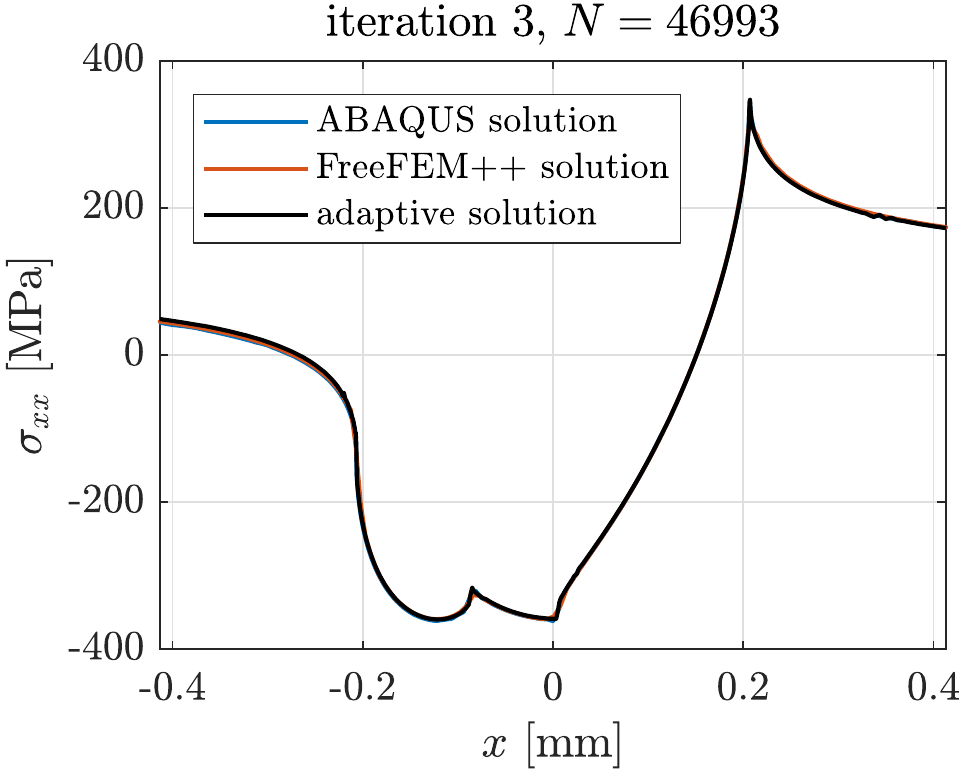}
  \caption{Surface traction $\sigma_{xx}$ under contact in four subsequent adaptive iterations during the solution of the fretting fatigue problem.}
  \label{fig:fwo-adapt-sxx}.
\end{figure}

The surface traction $\sigma_{xx}$ is of particular interest, as it is often used to
determine the location of crack initiation~\cite{pereira2016on}. Its graph over the course
of the adaptive iteration is shown in Fig.~\ref{fig:fwo-adapt-sxx}. In the initial solution,
only one node is placed under the contact surface; however, after the final iteration,
there are approximately 700 nodes under the contact surface and the results agree well with FEM
solutions.

\begin{figure}[h!]
  \centering
  \includegraphics[width=0.49\textwidth]{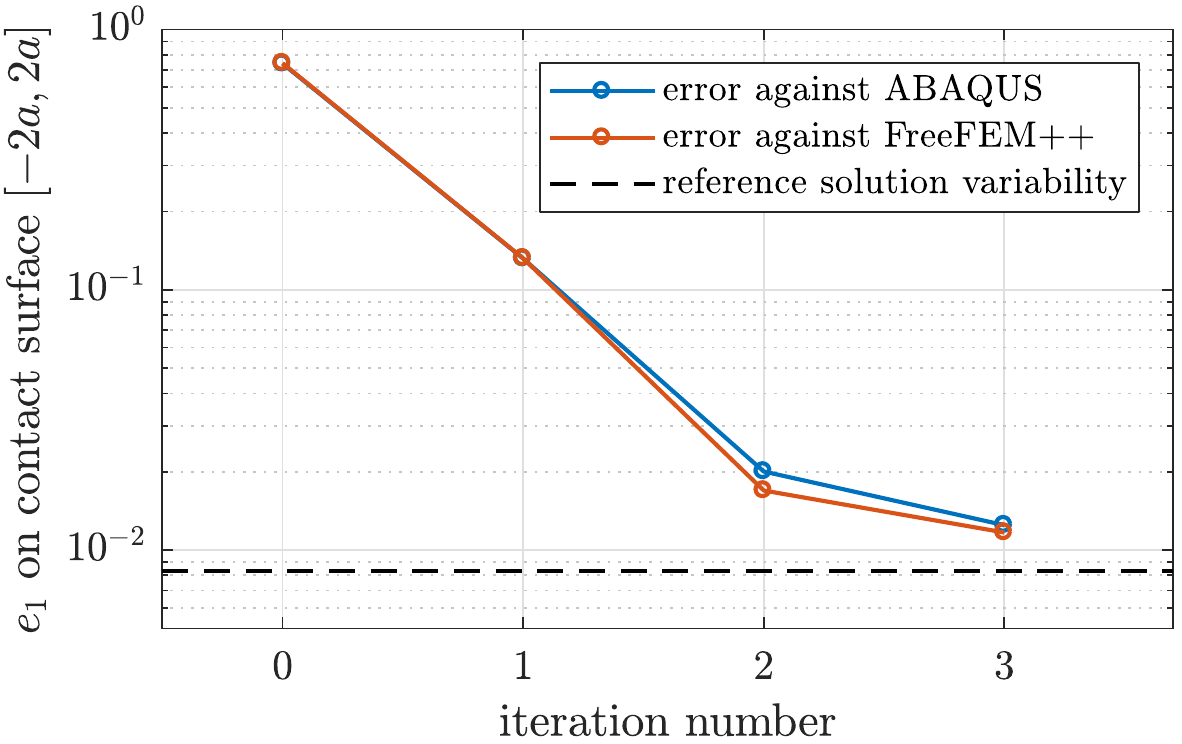}
  \caption{Error of surface traction using two reference solutions.}
  \label{fig:fwo-adapt-err}.
\end{figure}

Additionally, error of the adaptive solution is shown in Fig.~\ref{fig:fwo-adapt-err}.
The error was measured using $e_1$ (see Eq.~\ref{eq:err1}) on the contact
surface, specifically on the interval $[-2a, 2a]$, where $a$ is the contact half width.
As no analytical solution is known, the error was measured using reference FEM solutions.

The two reference solutions were computed to also asses the variability between solutions
obtained with established methods and implementations, and our adaptive implementation.
The error shown in Fig.~\ref{fig:fwo-adapt-err} decreases with each iteration, similarly to
other already presented cases, and achieves satisfactory accuracy, nearing the reference solution
variability.

Execution time was measured for both adaptive FEM and the proposed method.
All time measurements were done on a laptop computer in with
\texttt{Intel(R) Core(TM) i7-7700HQ CPU @ 2.80GHz}
processor and 16 GB DDR4 RAM. Code was compiled using \texttt{g++ (GCC) 8.1.1}
for Linux with \texttt{-std=c++11 -O3 -DNDEBUG} flags and no parallelisation.
The FEM solution was obtained using FreeFem++ version 3.61 compiled with
full optimisations on the same machine. The execution time varies
for both procedures with slight alterations of adaptivity parameters.
The observed execution times for solutions obtained with FreeFem++ ranged from
\unit[6]{s} to \unit[35]{s} and for the proposed adaptive procedure they
ranged from \unit[8]{s} to \unit[100]{s}.  For the specific solutions
used in Fig.~\ref{fig:fwo-adapt-err} the execution time were \unit[12.9]{s}
for FreeFem++ and \unit[38.6]{s} for the proposed procedure.
The reported times are total wall clock times which include generation of the initial discretization,
repeated refinement, matrix assembly, solution of the linear system and data export.
The majority of time is spent on solving the linear system, and in the meshless approach, on
computing the shape functions.

The execution times are reported only for reference and are not really comparable, since
adaptive procedures used are different and our implementation served to demonstrate a
proof-of-concept rather than being heavily focused on performance.

\subsection{Extension to 3-D problem}
\begin{figure}[h]
  \centering
  \includegraphics[width=0.4\linewidth]{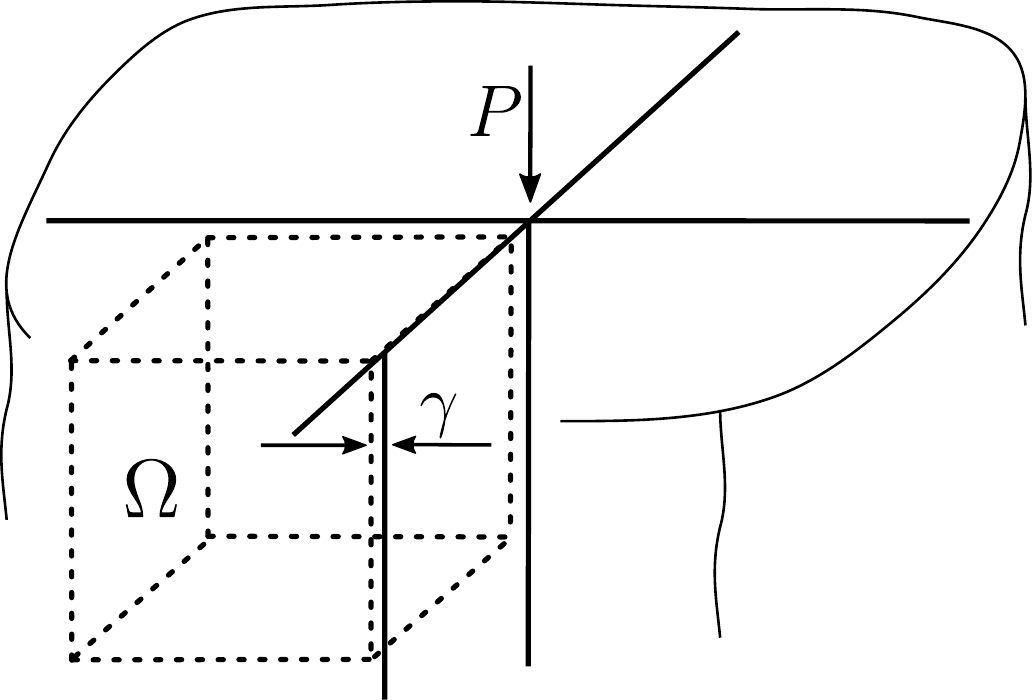}
  \caption{Boussinesq's problem and the computational domain $\Omega$.}
  \label{fig:3d-schema}
\end{figure}

\begin{figure}[h]
  \centering
  \includegraphics[width=0.4\linewidth]{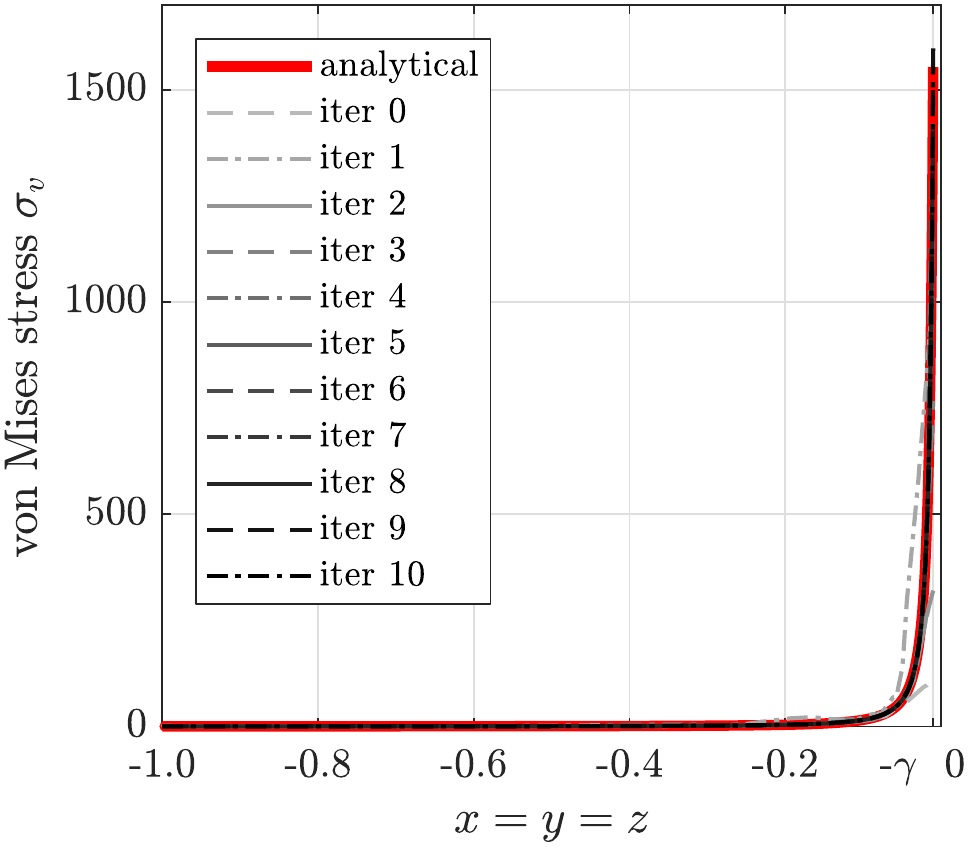}
  \includegraphics[width=0.4\linewidth]{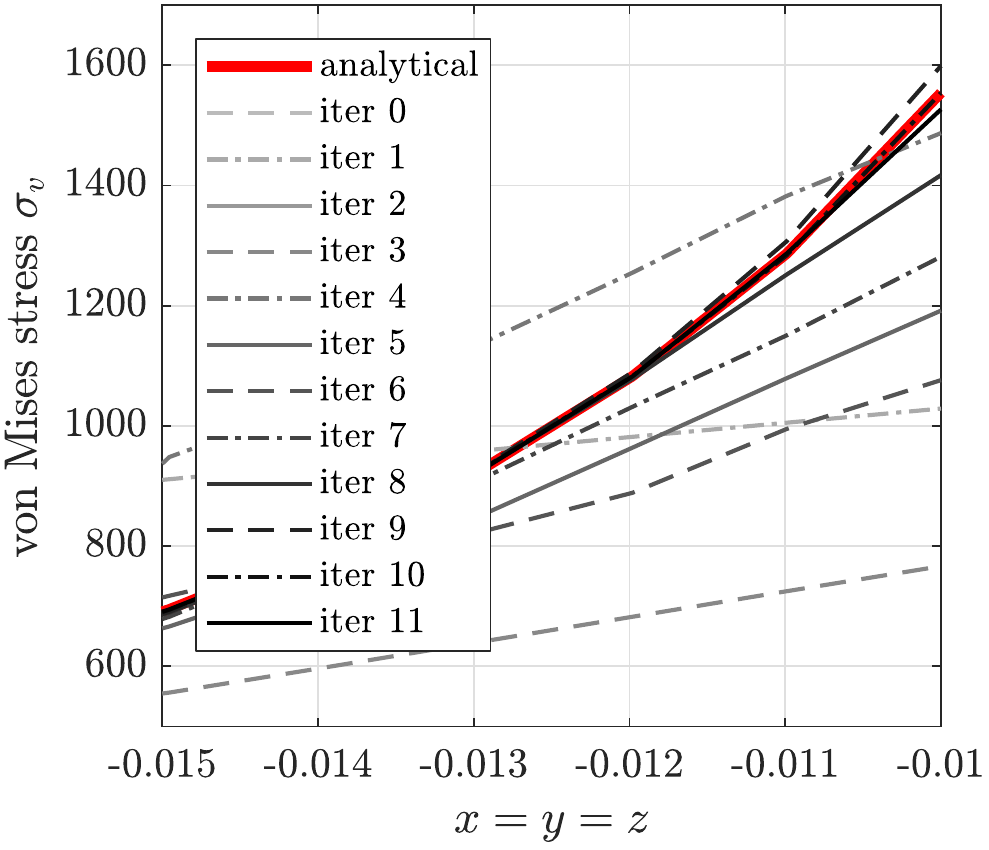}
  \caption{Stress profiles along the body diagonal from $[-1, -1, -1]$ to
  $[-\gamma, -\gamma, -\gamma]$ in adaptive iteration when solving the Boussinesq's problem.}
  \label{fig:3d-profile}
\end{figure}

In addition to the examples given in the previous sections, we present a fully
automatic adaptive solution of a 3-D problem. The Boussinesq's problem of the
concentrated normal traction acting on an isotropic half-space is considered as
described by Slaughter~\cite{Slaughter_2002}.
The closed form solution is given in cylindrical coordinates
$r$, $\theta$ and $z$ as
\begin{align}
  u_r &= \frac{Pr}{4\pi \mu} \left(\frac{z}{R^3} - \frac{1-2\nu}{R(z+R)} \right), \qquad
  u_\theta = 0, \qquad
  u_z = \frac{P}{4\pi \mu} \left(\frac{2(1-\nu)}{R} + \frac{z^2}{R^3} \right), \nonumber \\
  \sigma_{rr} &= \frac{P}{2\pi} \left( \frac{1-2\nu}{R(z+R)} - \frac{3r^2z}{R^5} \right), \qquad
  \sigma_{\theta\theta} = \frac{P(1-2\nu)}{2\pi} \left( \frac{z}{R^3} - \frac{1}{R(z+R)} \right),
    \label{eq:3d-problem} \\
  \sigma_{zz} &= -\frac{3Pz^3}{2 \pi R^5}, \qquad
  \sigma_{rz} = -\frac{3Prz^2}{2 \pi R^5}, \qquad
  \sigma_{r\theta} = 0, \qquad \sigma_{\theta z} = 0, \nonumber
\end{align}
where $P$ is the magnitude of the point force, $\nu$ is the Poisson's ratio, $\mu$ is the
second Lam\'e parameter and $R = \sqrt{r^2 + z^2}$ is the distance of given point from origin.
This solution has a singularity at the origin and similarly to previous sections, we will consider
a portion of a problem near the singularity. Numerically, we solve \eqref{eq:elasticity}
with essential boundary conditions given by~\eqref{eq:3d-problem} on a domain
\begin{equation}
\Omega = [-1, -\gamma] \times [-1, -\gamma] \times [-1, -\gamma]
\end{equation}
for $\gamma = 0.01$.
The domain $\Omega$ in an elastic half-space if shown in
Fig~\ref{fig:3d-schema}.
The similar stress profiles along the body diagonal of $\Omega$
are similar to the ones in Fig.~\ref{fig:anal-cross} and shown in
Fig.~\ref{fig:3d-profile}.

Although the solution is given in cylindrical coordinates, the discretisation
and enforcement of boundary conditions are done in Cartesian coordinates, with
the same  implementation of the method as used in 2-D is used in 3-D as well.
The domain was initially filled with constant density of $\drI = 0.05$,
amounting to $N = 7289$ nodes. For the given discretisation, RBF-FD
with 15 neighbouring nodes and 15 Gaussian RBFs was used. The shape parameter
$\sigma_b$ was varied proportionally to the internodal distance with base value
$\sigma_b = 100$.  For adaptive iteration
the straightforward generalisation of ad-hoc error indicator~\eqref{eq:err-indicator-std}
was used with $\eps = 5$, $\eta = 0$,
$\alpha = 3$ and $\beta = 1$. The reconstruction and evaluation of
the new density function $\widetilde{\dr}$ was done using $15$ closest
neighbours. Values of $P=-1$, $E=1$ and $\nu = 0.33$ were taken for
physical parameters of the problem.

\begin{figure}[h]
  \centering
  \includegraphics[width=0.5\linewidth]{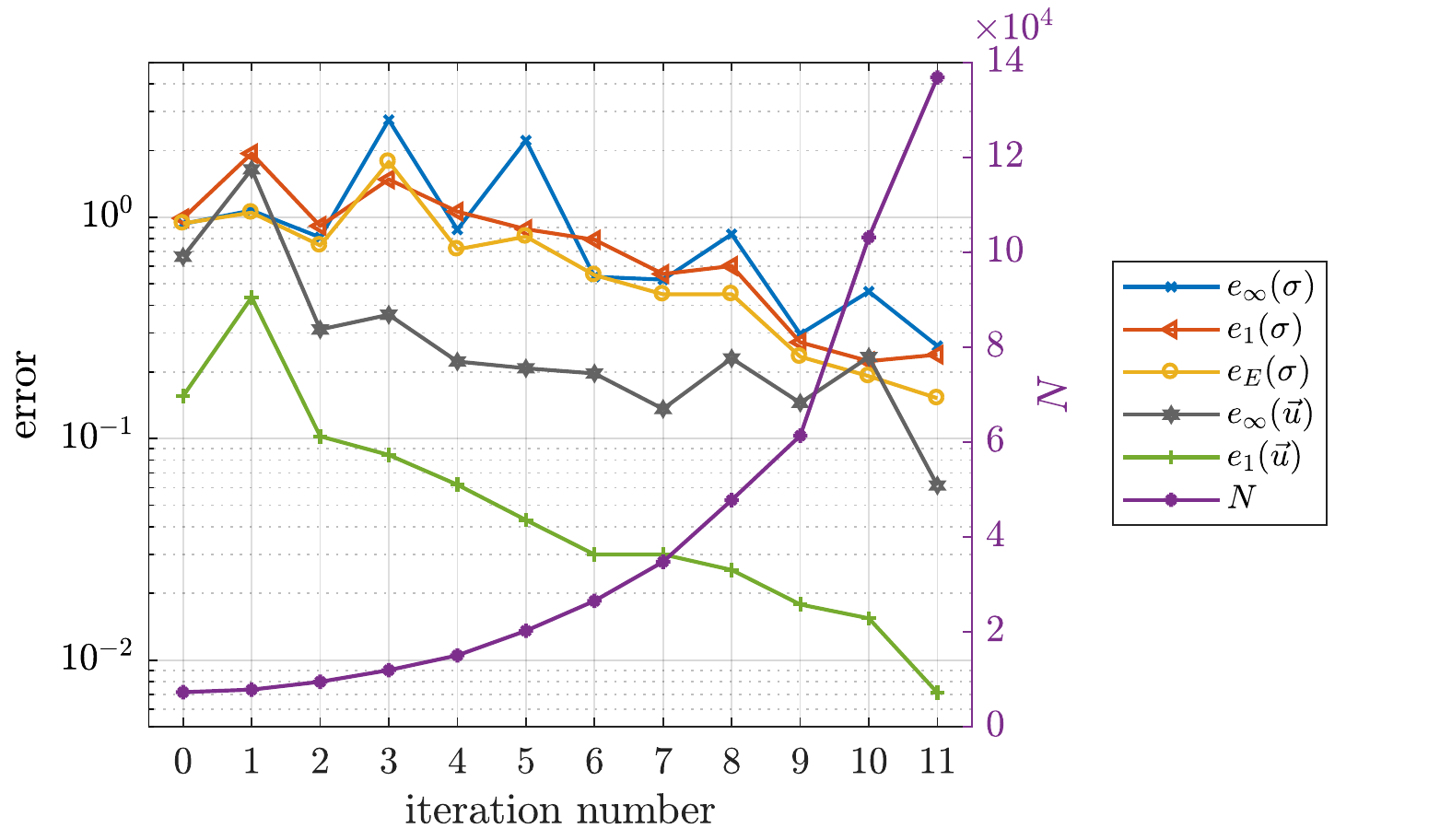}
  \caption{Error with respect to the number of nodes in the adaptive iteration
  of the solution of the Boussinesq's problem. The right axis shows the number
  of computational nodes.}
  \label{fig:error-3d}
\end{figure}

The error was measured using direct generalisations of formulas~(\ref{eq:errinf}--\ref{eq:energy-norm})
into 3-D. Since the closed form solution for displacement $\vec{u}$ is known in this case as well,
we will additionally show approximations of $L^1$ and $L^\infty$ errors for displacements as well,
which are computed analogously to~\eqref{eq:errinf} and~\eqref{eq:err1}
and will be denoted as $e_1(\vec{u})$ and $e_\infty(\vec{u})$, respectively.
Fig.~\ref{fig:error-3d} shows the computed errors over the course of the adaptive iteration.
Additionally, von Mises stress profiles along the body diagonal from $[-1, -1, -1]$ to
$[-\gamma, -\gamma, -\gamma]$ are shown in Fig.~\ref{fig:3d-profile}.

The behaviour of the method is similar to the two-dimensional compressed disk case.
Initially, the error rises until the area around the contact is covered densely enough.
The displacement errors are smaller than stress errors, as is customary in solid mechanics.
Right panel in Fig.~\ref{fig:3d-profile} shows more precisely how the peak stress is estimated.

\begin{figure}[h]
  \centering
  \includegraphics[width=0.45\linewidth]{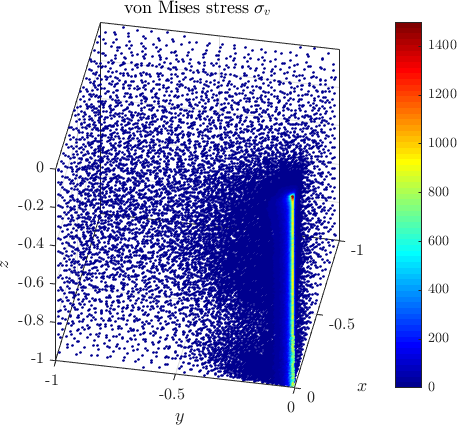}
  \includegraphics[width=0.4\linewidth]{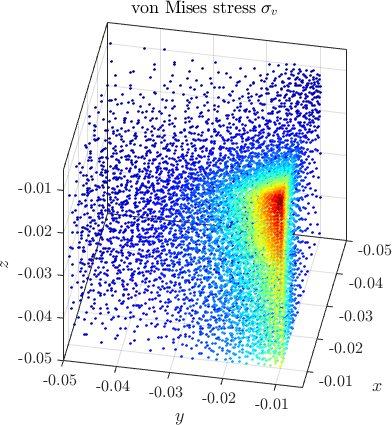}
  \caption{The obtained solution in the final iteration
    with an enlarged portion around the contact area. Both solutions are plotted only
    in the computation nodes to also show the final nodal distribution. Nodes are coloured
    proportional to computed values of von Mises stress.}
  \label{fig:sol-3d}
\end{figure}

Fig.~\ref{fig:sol-3d} shows the solution in the final iteration of the adaptive procedure.
The corner $[-\gamma, -\gamma, -\gamma]$ is shown more closely in the right figure.

The execution time of the entire adaptive procedure to obtain the solution shown in Fig.~\ref{fig:sol-3d}
was \unit[125]{s} with \unit[1]{s} spent on the initial iteration with 7289 DOFs and
\unit[48]{s} spent on the last iteration with 138\,632 DOFs. In the last iteration, the
internodal distances between the densest parts are about 176 times smaller than in
the coarsest parts.

\section{Conclusions}
\label{sec:conclusions}
Adaptive refinement is a crucial part of a numerical solution procedure for any
problem that exhibits extensive differences in intensity within the problem
domain. In this paper, a modular adaptive refinement algorithm for solving
elliptic boundary value problems was proposed. The refinement procedure consists
of a Poisson Disc Sampling based nodal positioning algorithm, an ad-hoc error
indicator, a RBF-FD meshless solution of the governing problem, the refinement
strategy, and the reconstruction of the density function. The
free parameters of the algorithm were determined by numerical
experiments, and reasonable defaults along with acceptable intervals were
proposed. Finally, it was clearly demonstrated that the proposed algorithm
provides a stable solution for all considered cases without any hyper-tuning.
This includes the solution of the fretting fatigue case, where complications
arise due to the presence of the stick-slip zone at the contact interface, which
causes extreme peaks in surface stress making it very difficult to solve
numerically in a reasonable time without refinement. The proposed approach
successfully adaptively refined initially uniform discretisation, captured both
peaks, and numerically solved the problem in only four iterations, with
excellent agreement towards the reference solution provided by commercial
software using manual refinement.

All elements of the algorithm have no intrinsic dependency on the dimensionality
of the space and generalisation to higher dimensions is possible. In
the final example, the generality of the proposed method is demonstrated in the
solution of 3-D Boussinesq's problem of the concentrated normal traction acting
on an isotropic half-space.

Although all discussions in this paper were limited to the consideration of the
Navier-Cauchy partial differential equation and contact problems, the presented
methodology can also be applied in other areas where numerical discretisation is
required to solve the governing problem. The only notable required change would
be in the second step, i.e.\ the solution of the problem, while the other steps
would remain the same.

The examples in this paper were computed using the in-house open source Medusa
library~\cite{medusa} and the algorithms presented in this paper. All research
data is freely available on a git repository~\cite{git}.

Future work will be focused on implementation of more error indicators, more
efficient implementation, and extension to different problems.

\section{Acknowledgements}
Authors would like to acknowledge the financial support of the Research
Foundation Flanders (FWO), The Luxembourg National Research Fund (FNR) and
Slovenian Research Agency (ARRS) in the framework of the FWO Lead Agency
project: G018916N Multi-analysis of fretting fatigue using physical and virtual
experiments, and the ARRS research core funding No.\ P2-0095.

\bibliography{references}%

\end{document}